\documentclass[12pt,twoside]{article}
\usepackage{color}
\usepackage{indentfirst}
 \usepackage{amsmath,amssymb,amscd}
\usepackage{cite}
 \newcommand{\Ol}{{\Omega_L}}

\newtheorem{lemma}{\bf Lemma}[section]
\newtheorem{theorem}{\bf Theorem}
\newtheorem{proposition}[lemma]{\bf Proposition}
\newtheorem{corollary}[lemma]{\bf Corollary}
\newtheorem{defi}[lemma]{\bf Definition}
 \newtheorem{rem}[lemma]{Remark}

       \newtheorem{remark}[lemma]{\bf Remark}

\newcommand{\btt}{\begin{theorem}}
\newcommand{\bl}{\begin{lemma}}
\newcommand{\el}{\end{lemma}}
\newcommand{\ett}{\end{theorem}}

\newcommand{\no}{\nonumber\\}
\newcommand{\rr}{\rightarrow}
\newcommand{\pa}{\partial}

\newcommand{\bi}{\bibitem}

\newcommand{\vp}{\varphi}
\newcommand{\dis}{\displaystyle}

\newcommand{\ga}{\gamma}
\newcommand{\te}{\theta}

\newcommand{\al}{\alpha}
\newcommand{\de}{\delta}
\newcommand{\ve}{\varepsilon}

\newcommand{\aaa}{\,\,\,\,\,\,}

 \newcommand{\lb}{\label}
\newcommand{\ol}{\overline}

\newcommand{\bn}{\begin{eqnarray}}
\newcommand{\en}{\end{eqnarray}}
\newcommand{\bnn}{\begin{eqnarray*}}
\newcommand{\enn}{\end{eqnarray*}}

 \let\dis=\displaystyle
 
 \let\al=\alpha
 \let\be=\beta
 \let\ga=\gamma
 \let\de=\delta
 \let\eps=\varepsilon

 \let\pa=\partial

 \let\si=\sigma
 \let\Om=\Omega

\def\endsquare{\ifmmode \kern 0.5em\square \else \discretionary
 {}{\hbox to\hsize {\hfill \square }\null } {\kern 0.5em\square }\fi }

\def\r  { \rho_\varepsilon }

\def\uu {{u_\varepsilon}}

 \def\be#1\ee{\begin{equation}#1\end{equation}}
 \def\bay#1\eay{\!\!\!\left\{\!\!\begin{array}{l}%
      #1\displaystyle\end{array}\right.}
 \def\bln#1\eln{\begin{aligned}#1\end{aligned}}

 \def\bma#1\ema{{\allowdisplaybreaks\begin{align}#1\end{align}}}
 \def\nnm{\notag}
 \def\bgr#1\egr{{\allowdisplaybreaks\begin{gather}#1\end{gather}}}

       \def\ef#1{(\ref{#1})}
       \def\qef#1{$(\ref{#1})$}

       \def\oplem#1{\begin{lemma}\, {\rm #1}\, \it }
       \def\cllem{\end{lemma}\rm \par }

       \def\opthm#1{\begin{theorem}\, {\rm #1}\, \it }
       \def\clthm{\end{theorem}\rm \par }

       \def\oppro#1{\begin{proposition}\, {\rm #1}\, \it }
       \def\clpro{\end{proposition}\rm \par }

        \def\opcor{\begin{corollary} }
       \def\clpcor{\end{corollary}\rm \par }

       \def\opdef{\begin{definition}}
       \def\cldef{\end{definition}\rm\par }

       \def\oprem{\begin{remark}\it}
       \def\clrem{\end{remark}\rm \par}

       \def\demo{{\bf Proof: \ }\rm }
       \def\enddemo{\null\hfill\framebox{}\par\bigskip}

\def\n{\rho}
\begin{document}
\title{Vanishing of Vacuum States and Blow-up Phenomena
   of the Compressible Navier-Stokes Equations}
 \author{\textbf{Hai-Liang LI$^{1,3}$,\quad
Jing LI$^{2,3}$,\quad
 Zhouping XIN$^{3,1}$}\\[2mm]
 {\small\it\it $^{1}$Department of Mathematics,
       Capital Normal University}\\
 {\small\it Beijing, P. R. China,
 email: hailiang$_-$li@mail.cnu.edu.cn}\\
{\small\it\it $^{2}$Institute of Applied Mathematics,
  AMSS, Academia Sinica}\\
 {\small\it Beijing 100080,
 P. R. China ,
 email: ajingli@gmail.com}\\
{\small\it\it $^{3}$Institute of Mathematical Science,
 The Chinese University of Hong Kong}\\
 {\small\it Shatin, Hong Kong,
 email: zpxin@ims.cuhk.edu.hk}
}

\pagestyle{myheadings} \markboth{Vanishing of vacuum states and
blow-up phenomena} {H.-L. Li, J. Li, $\&$ Z. Xin}
\date{}
\maketitle

\begin{abstract}

The Navier-Stokes systems for compressible fluids with
density-dependent viscosities are considered in the present paper.
 These equations, in particular, include the ones which are
rigorously derived recently  as the Saint-Venant system for the
motion of shallow water, from the Navier-Stokes system for
incompressible flows with a moving free
surface~\cite{GerbeauPerthame01}. These compressible systems are
degenerate when vacuum state appears. We study
initial-boundary-value problems for such systems for both bounded
spatial domains or periodic domains. The dynamics of weak
solutions and vacuum states are investigated rigorously.

First, it is proved that the entropy weak solutions for general
large initial data satisfying finite initial entropy exist globally
in time. Next, for more regular initial data,   there is a global
entropy weak solution which is unique and regular with well-defined
velocity field for short time, and the interface of initial vacuum
propagates along particle path during this time period. Then, it is
shown that for any global  entropy weak solution, any (possibly
existing) vacuum state must vanish within finite time. The velocity
(even if regular enough and well-defined) blows up in finite time as
the vacuum states vanish. Furthermore, after the vanishing of vacuum
states, the global entropy weak solution becomes a strong solution
and tends to the non-vacuum equilibrium state exponentially in time.
\end{abstract}



\section{Introduction}\setcounter{equation}{0}

The compressible isentropic  Navier-Stokes equations, which are the
basic models describing the evolution of a viscous compressible
fluid,
  read as follows \bn
\lb{ii1}\begin{cases}\rho_t+div(\rho u)=0,\\
 (\rho u)_t+div(\rho u\otimes u)- 2div(\mu  D( u)
 )-\nabla(\xi
div u)+\nabla p(\rho)=0 ,\end{cases}\en where $x\in\Om\subset R^N,
t\in (0,T),$ $D(u)=(\nabla u+(\nabla u)^{\rm T})/2,$ and $
p(\rho)=a\rho^\ga,a>0,\ga\ge 1,$ the viscosity coefficients
$\mu,\xi$ are assumed to satisfy $\mu\ge 0$ and $\xi+2\mu/N\ge 0.$

If    $\mu$ and $\xi$ are both constants, there is huge literature
on the studies of the global existence and behavior of solutions to
\ef{ii1}. For instance, the one-dimensional (1D) problems were
addressed by Kazhikhov et al \cite{ks1} for sufficiently smooth
data, and by Serre \cite{se1,se2} and Hoff \cite{ho1} for
discontinuous initial data, where the data are uniformly away from
the vacuum; the multidimensional problems (\ref{ii1}) were
investigated by Matsumura et al \cite{ma2,ma1,ma3}, who proved
global existence of smooth solutions for data close to a non-vacuum
equilibrium, and later by Hoff  for discontinuous initial data
 \cite{ho3}, and more recently, by Danchin \cite{da1},
who obtained existence and uniqueness of global solutions in a
functional space  invariant by the natural scaling of the
associated equations; and for the existence of solutions for
arbitrary data(which may include vacuum states), Lions
\cite{li3,li4,Lions98} (see also Feireisl et al \cite{fe1})
obtained global existence of weak solutions - defined as solutions
with finite energy - when the exponent $\ga$ is suitably large,
where the only restriction on initial data is that the initial
energy is finite, so that the density is allowed to vanish.

Despite the important progress, the regularity, uniqueness and
behavior of these weak solutions remain largely open. As
emphasized in many papers related to compressible fluid dynamics
\cite{chok1,cho2,ho1,HoffSerre91,HoffSmolller01,ks1,
li04,sa1,se1,so2,vz1,Xin98a}, the possible appearance of vacuum is
one of the major difficulties when trying to prove global
existence and strong regularity results. Hoff and Smoller
\cite{HoffSmolller01} proved that weak solutions of the
compressible  Navier-Stokes equations \ef{ii1} in one space
dimension do not exhibit vacuum states in a finite time provided
that no vacuum is present initially under fairly general
conditions on the data.  Such a  result was  extended to the
spherically symmetric case in \cite{XinYuan05} recently. On the
other hand, the results of Xin \cite{Xin98a} showed that there is
no global smooth solution  to Cauchy problem for (\ref{ii1}) with
a nontrivial compactly supported initial density, which gives
results for finite time blow-up in the presence of vacuum. It is
also proved in  \cite{li04} that for bounded domain even one point
initial vacuum shall cause  the global strong solutions to 1D
\ef{ii1} to blow up as time goes to infinity provided the initial
data satisfies some compatibility conditions.

 The independence of viscosities on density  makes it possible to
trace either the particle path or the trajectory of vacuum state.
This, however, leads to the failure of continuous dependence of weak
solutions containing vacuum state on initial data
\cite{HoffSerre91}. For the case that the density changes
continuously across the interfaces separating the gas  and vacuum,
the global existence and uniqueness of  weak solutions was obtained
in \cite{LuoXinYang00}, where the authors obtained that velocity is
smooth enough up to the interfaces which are particle paths
separating the gas from the vacuum and that the support of gas
density expands outside as well as interface connecting gas and
vacuum moves at an algebraic rate.
\par

 Thus, viscous compressible fluids near vacuum should be
 better
modeled by the compressible Navier-Stokes equations with
density-dependent viscosities, as was derived in the fluid-dynamical
approximation of Boltzmann equation for dilute gases. Further, as
was first pointed out and investigated by Liu-Xin-Yang in
\cite{LiuXinYang98}, in the  derivation of  the compressible
Navier-Stokes equations from the Boltzmann equation by  the
Chapman-Enskog expansions, the viscosity  depends on the
temperature,which is translated into the dependence of the viscosity
on the density for isentropic flows. Moreover, it should be
emphasized that a one-dimensional  compressible flow model, called
the viscous Saint-Venant system for laminar shallow water, derived
rigorously from incompressible Navier-Stokes system   with a moving
free surface  by Gerbeau-Perthame recently
in~\cite{GerbeauPerthame01}, has the form: \be
\lb{ii2}\begin{cases}\rho_t+ (\rho u)_x=0,\\
 (\rho u)_t+ (\rho u^2)_x -a(\n u_x)_x + (\rho^2)_x=0 ,\end{cases}
 \ee
with the viscosity coefficients given by $\mu(\n)=\xi(\n)=a\n/3$
for a given positive constant $a.$ Indeed, such models appear
naturally and often in geophysical flows
\cite{BreschDesjardins03,BreschDesjardins05,BreschDesjardinsLin05}.

In the  case of one-dimensional problem with  $\mu=\xi=a\n^\al$ for
some positive constants $a$ and $\al,$   the well-posedness of the
Cauchy problem has been studied by many authors for not only
initially compact-supported density but initially non-vacuum
density.  Indeed, for the initially compact-supported density case,
the local (in time) well-posedness of weak solutions to this problem
was first established by Liu-Xin-Yang in \cite{LiuXinYang98}, where
the initial density was assumed to be connected to vacuum with
discontinuities. This property, as shown in \cite{LiuXinYang98}, can
be maintained for some finite time. And the global existence of weak
solutions, together  with $\al\in(0,1)$ and the density function
connecting to vacuum with discontinuities, was considered by many
authors, see \cite{OkadaMatsusuMakino02,yangyao01,JiangXinZhang06}
and the references therein. It is noticed that the above analysis is
based on the uniform positive lower bound of the density with
respect to the construction of the approximate solutions.
On the other hand,  if the density function connects to vacuum
continuously, there is no positive lower bound for the density and
the viscosity coefficient vanishes at vacuum. This degeneracy in
the viscosity coefficient gives rise to new difficulties in
analysis because of the less regularizing effect on the solutions.
Yang et al \cite{Yangzhao02} first obtained a local existence
result for this case under the free boundary condition with
$\al>1/2.$ The authors in \cite{Yangzhu02,vong,FangZhang04}
obtained the global existence of weak solutions with
$\al\in(0,1/2). $

However, almost all above results   concern mainly with free
boundary problems, and for the global existence results, the choices
of viscosity do  not fit the important physical model, the shallow
water equation \ef{ii2} with $\mu(\n)=\n$ (namely $\al=1$). For the
constant viscosity case, one has known not only that the vacuum
state will not develop later on time if there is no vacuum state
initially~\cite{HoffSmolller01}, but also that the separate two
initial vacuum states shall not meet together in a finite
time~\cite{XinYuan05},  and  that one   point initial vacuum causes
strong solutions  to blow up~\cite{li04} at infinity as well.
However, little is known  on   the dynamics of the vacuum states  of
weak solutions to the compressible Navier-Stokes equations \ef{ii1}
with density-dependent viscosity on bounded domain. And in
particular,
 it is not clear yet how the vacuum states evolve  with
respect to time and whether the initial vacuum states shall exist
all the time or not for weak solutions. The study of these important
dynamical problems about vacuum states is rather difficult because
the nonlinear diffusion is degenerate as vacuum appears, which is
quite different from the case of constant viscosity. This causes the
loss of information about the velocity and makes it difficult to
trace the evolution of vacuum states in general.  Across the
interface (or vacuum boundary), it is usually difficult to obtain
enough information about velocity even if considering specific cases
such as point vacuum or continuous vacuum of one piece. It is
important to get enough information about  the velocity since the
flow particles transport usually along particle path, and all
interfaces of vacuum, such as free
boundaries~\cite{LuoXinYang00,JiangXinZhang06} which can be observed
and dealt with, also move along the trajectories determined by
velocity field.

For the multidimensional case, Vaigant et al \cite{vk95} first
proved that for the 2D case and for the case $\mu$ is a constant and
$\xi(\n)=a\n^\beta,$ with $a>0,\beta> 3, $   \ef{ii1} with periodic
boundary condition has a unique strong and classical solution with
density away from vacuum. More recently, Bresch  and Desjardins
\cite{BreschDesjardinsLin05,BreschDesjardins03,
BreschDesjardins03b,BreschDesjardins05} (see also
\cite{MelletVasseur05}) have made important progress. Under the
condition that $\xi(\n)=2(\mu'(\n)\n-\mu(\n)),$ they establish a new
Bresch-Desjardins (BD) entropy  inequality which can not only be
applied to the vacuum case but also used to get the compactness
results for \ef{ii1} which extended the compactness results due to
Lions \cite{li3,li4,Lions98} to the case $\ga\ge 1$.
 On the other hand, the constructions of the approximation solutions
does not seem routine in the general case of appearance of vacuum.
However, it should be noted that recently Bresch et al
\cite{bd07,bdg07,BreschDesjardins05} have made significant progress
on the construction of approximate solutions and existence of global
weak solutions to the multi-dimensional compressible Navier-Stokes
equations and the 2D shallow water model in the case that there is
either a drag friction (for barotropic compressible flows) or a cold
pressure (for viscous and heat conducting flows). In the case that
there is neither drag friction nor cold pressure included,
Guo-Jiu-Xin~\cite{GuoJiuXin2007} recently have shown how to
construct approximate smooth solutions and obtain the global
existence of weak solutions via the BD entropy to the (2D and 3D)
barotropic compressible Navier-Stokes for the spherical symmetric
initial data. There are also other recent interesting applications
of the BD entropy to one-dimensional compressible Navier-Stokes
equations with degenerate viscosities, for instance, on the global
existence and long time behavior of weak solutions for free boundary
problem~\cite{GuoZhu2007,GuoJiangXie2007}, or the existence and
uniqueness of global strong solution away from vacuum in real
line~\cite{mv07}.
\par

We study  mainly the initial-boundary-value problem (IBVP) for
\ef{ii1}, where $\mu=\n^\al$ with $\al> 1/2,$ on spatial
one-dimensional bounded spatial domains or periodic domains. This
contains the physical important model for shallow water equations
\ef{ii2}. The choice of $\al> 1/2$ is necessary in order to consider
the dynamics of vacuum states since it allows the existence of
initial vacuum in Eulerian coordinates as one can see later.

We first establish the global existence of entropy weak solutions
for the compressible Navier-Stokes equations \ef{ii1}, with pressure
$p= \n^\ga$ and $\ga\ge \alpha/2$, for general initial data with
finite entropy and vacuum. The key in our analysis is the
construction non-vacuum approximate solutions so that we can make
use of the stability analysis in \cite{MelletVasseur05}, where the
Bresch-Desjardins (BD) entropy inequality (see
\cite{BreschDesjardinsLin05,BreschDesjardins03,
BreschDesjardins03b,BreschDesjardins05}) was used to obtain  the
compactness results.
Our construction of the approximate solutions is strongly motivated
by the previous work of Jiang-Xin-Zhang~\cite{JiangXinZhang06} about
the existence of global weak solutions to one-dimensional
compressible Navier-Stokes with free boundary as vacuum interface.
In general, it seems rather difficult to investigate the dynamics of
vacuum states due to the degeneracy of nonlinear diffusion and the
density function connecting to vacuum continuously. Therefore, we
further consider the cases of more regular initial data containing
point vacuum or continuous vacuum of one piece, and we show that
there is a global
 entropy weak solution which is unique and regular with
well-defined velocity field at least for short time, and the vacuum
states remain for the short time. Then, we use some ideas due to
\cite{LiXinHuang05,li04,LiXin04} to prove that any possible vacuum
state in such global weak solutions which satisfy the BD entropy
must vanish within finite time. This shows that such short time
structure and vacuum states of weak solutions can not be maintained
all the time. And as the vacuum states vanish, the spatial
derivative of velocity (if it exists) has to blow up even if the
velocity is regular enough and well-defined before. After the
vanishing of vacuum states, we can redefine the velocity field and
recover the nonlinear diffusion term in terms of density and
velocity. In addition, the global entropy weak solution is shown to
become a strong solution and tends to the non-vacuum equilibrium
state exponentially in time. This phenomena, applied to the
compressible shallow water equations \ef{ii2},  seems to be never
observed for the compressible Navier-Stokes equations before.
\bigskip

The rest  of the paper is as follows. In section~\ref{main}, the
main results about the vanishing of vacuum states and blow-up
phenomena of global entropy weak solutions for the compressible
Navier-Stokes equations are stated. The global existence of
entropy weak solutions for general large initial data with vacuum
states allowed is proven in section~\ref{Global-Theory}. The short
time structure of global entropy weak solution with initial one
point vacuum state or initial continuous vacuum states of one
piece are investigated in section~\ref{Dynamics}. In
section~\ref{Vanishing}, we show the vanishing of vacuum states
and blow-up phenomena of any global entropy weak solution within
finite time and analyze the regularity and large time asymptotic
behavior of global entropy weak solutions after the vanishing of
vacuum states.

\section{Main results}\lb{main}\setcounter{equation}{0}
We consider the initial-boundary-value-problem (IBVP) for the 1D
compressible Navier-Stokes equations with density-dependent
viscosity \bgr
 \n_t + (\n u )_x = 0, \label{l0a}
\\
 (\n u )_t + (\n u^2 + p(\n))_x
            -(\mu(\n)u_x)_x=0,
  \label{l0b}
 \egr
with $\n\ge0$ the density, $\n{u}$ the momentum. The pressure and
viscosity are assumed to have the form: \bnn
 p(\n)=a_1\n^\gamma,\ \mu(\n)=a_2\n^\alpha
  \enn
where $\gamma \ge 1,\ a_1>0,\ a_2>0,$ and $ \alpha> 1/2$ are
constants, and for simplicity we set $a_1=a_2=1$.

The initial data is given for the density $\n$ and the momentum
$\n{u}$ \be \n(x,0)=\n_0(x)\ge 0,\quad
 \n u(x,0)= m_0 (x),  \quad x\in \Om,  \label{l0c}
\ee
where the domain $\Om$ is chosen as unit interval denoting the
spatial domain $(0,1)$ or periodic domain with period length one,
and throughout the present paper the initial data is assumed to
satisfy
\be
\begin{cases}
\,  \n_0\ge 0 \mbox{ a.e. in }\Om,\quad
   \n_0\in L^1(\Om),\quad
  (\rho_0^{\al-1/2})_x\in L^2(\Om),\\[2mm]
\, m_0=0,\  \mbox{ a.e. on }\{x\in \Om\,|\,\n_0(x)=0\},\
  \mbox{$\frac{|m_0|^2}{\n_0}$}\in L^{1}(\Om).
\end{cases}\lb{d8}
\ee \oprem{} Note here that the condition \qef{d8} implies \be
\n_0\in L^\infty(\Om),\quad  \n_0\log_+\n_0\in L^1(\Om). \lb{d8zz}
\ee
 It should be clear that a large class of initial data satisfy the
 conditions in \qef{d8}. In particular, the assumptions \qef{d8} are satisfied for
following initial data
 \[
\n_0(x)=(|x-x_0|^2)^{1/(2\alpha-1)},\quad m_0(x)=0,\quad x\in\Om.
 \]
Without the loss of generality, the total initial mass is
renormalized to be one throughout the present paper, i.e., \be
 \int_\Om\n_0(x)dx=1.   \nnm\label{l0dIDb}
\ee
\clrem The boundary conditions are one of the boundary conditions
of Dirichlet type and periodic type for Eqs.~\ef{l0a}-\ef{l0b}
imposed as
 \begin{description}
\item{(1).}  Dirichlet case:
 \be
\n{u}(0,t)=\n{u}(1,t)=0,\quad t\ge 0, \label{d6}
\ee

\item{(2).} periodic case:
\be
\n,\ u \mbox{ are periodic in } x \mbox{ of
period }1,\lb{d7}
\ee
where we consider Eqs.~\ef{l0a}-\ef{l0b} on
$\mathbb{R}\times(0,\infty)$.
\end{description}
\oprem{}
For the case of Dirichlet boundary \qef{d6}, the boundary is given
by the physical observable momentum instead of velocity. It is
natural to employ such boundary condition as considering the
dynamics of (global in time) weak solutions to  the IBVP  for the
compressible Navier-Stokes equations with possible vacuum states
included since it is usually the case that, at vacuum states, the
momentum is zero and is observed and controllable, but almost
nothing is known yet for the velocity for weak solutions. Note
here that for any (weak) solution away from vacuum at the
boundary, the boundary condition \qef{d6}   reduces to
$u(0,t)=u(1,t)=0$, $t\ge 0$. \clrem In order to define the weak
solutions to the IBVP for the compressible Navier-Stokes
equations~\ef{l0a}-\ef{l0b} with initial data \ef{l0c} and
boundary condition \ef{d6} or \ef{d7}, we define the set of test
functions as follows, \bnn \Psi\triangleq\begin{cases}C_0^\infty(
\ol\Om\times[0,T)) &\mbox {
for the Dirichlet case }\ef{d6},\\
   C_{\rm per} ^\infty(\mathbb{R}\times[0,T)) &\mbox { for
the periodic case }\ef{d7},\end{cases}\enn and \bnn
\Phi\triangleq\begin{cases}C_0^\infty( \Om\times[0,T)) &\mbox {
for the Dirichlet case }\ef{d6},\\
   C_{\rm per} ^\infty(\mathbb{R}\times[0,T)) &\mbox { for
the periodic case  }\ef{d7},\end{cases}\enn  with $C_{\rm per}
^\infty(\mathbb{R}\times[0,T))$  defined by\bnn C_{\rm per}
^\infty(\mathbb{R}\times[0,T))= \left\{\vp\in
C^\infty(\mathbb{R}\times[0,T)) | \,\,\vp  \mbox{  is periodic in
$x$ of period  1}\right\}.
 \enn
We define the  weak solutions to  the IBVP for the compressible
Navier-Stokes Equations~\ef{l0a}-\ef{l0b} as follows.
\begin{defi}[global weak solutions]\label{definition}
For any $T>0,$ $(\rho,u)$ is said to be a  weak solution to
Eqs.~\qef{l0a}-\qef{l0b} with  initial data \qef{l0c} and boundary
value \qef{d6} or \qef{d7}  in $\Om\times(0,T)$, if  \bn
\begin{cases}
 0\le \rho\in L^\infty(0,T;L^1(\Om)\cap L^\gamma(\Om) ),\quad
 (\rho^{\al-1/2})_x\in L^\infty(0,T;L^2(\Om)), \\[2mm]
  \sqrt {\n }u\in L^\infty(0,T;L^2(\Om)),\quad
 \rho^{\al } u_x\in L^2(0,T;W^{-1,1}_{\rm loc}(\Om)),
\end{cases} \lb{s1}
\en  and $(\n,u)$ satisfies
 \bn
 \int_\Om \n_0\psi(x,0) dx +\int_0^T\int_\Om \rho\psi_tdxdt
 +\int_0^T\int_\Om\sqrt{\n}\sqrt{\n}u \psi_xdxdt=0        \lb{d1}
\en for any $\psi\in \Psi $, and
 \bn
 \lefteqn{\int_\Om m_0\vp(x,0)dx
 +\int_0^T\int_\Om \sqrt{\rho}(\sqrt{\rho}u) \vp_t
 dxdt}  \no&& +\int_0^T\int_\Om
\left( (\sqrt{\rho}u)^2+\rho^\ga\right)\vp_xdxdt- \langle\n^\al
u_x,\vp_x\rangle =0
 \lb{d2}
\en for all $\vp\in\Phi$.  The nonlinear diffusion term $
 \rho^{\al } u_x $ is defined as \bn
  \langle\n^\al
u_x,\vp \rangle
 &=&
 -\int_0^T\int\rho^{\al-1/2}\sqrt{\rho}u\vp_{x }dxdt\no&&
 - \frac{2\al}{2\al-1}
   \int_0^T\int(\rho^{\al-1/2})_x\sqrt{\rho}u\vp dxdt\lb{d3}
\en
 for any $\vp\in \Phi$, where  $\n\in L^\infty(\Om\times(0,T))$
due to \qef{s1}. Moreover, for the spatial periodic case \qef{d7},
$(\rho,\sqrt{\n}u)$ is also periodic.
\end{defi}

\begin{remark} For the Dirichlet case, \qef{d1}, together with
the fact $\n\in L^\infty(\Om\times(0,T))$ due to \qef{s1}, implies
that $(\n,u)$ satisfies the Eq.~\qef{l0a} in the sense of
distribution and justifies the boundary condition \qef{d6} in the
sense that, for any time interval $I\subset [0,T],$
$$\ve^{-1}\int_I\int_0^\ve\sqrt{\n}\sqrt{\n}u(x,s)dxds\to 0,\quad
\ve^{-1}\int_I\int_{1-\ve}^\ve\sqrt{\n}\sqrt{\n}u(x,s)dxds\to 0$$ 
 as $\ve\to 0^+$. If further
$\sqrt{\n}\sqrt{\n}u\in L^p(0,T;W^{1,q} (\Om))$ for some $p\ge
1,q\ge 1,$ then \qef{d1} yields that $\sqrt{\n}\sqrt{\n}u\in
L^p(0,T;W^{1,q}_{0}(\Om )),$ that is, the Dirichlet boundary
condition \qef{d6} is satisfied in the sense of trace.
\end{remark}

\begin{defi}[global  entropy weak solutions]\label{definition1}
Let  $(\rho,u)$ be a global weak solution (in the sense of
Definition \ref{definition}) to \qef{l0a}-\qef{l0b} with  initial
data \qef{l0c} and boundary value \qef{d6} or \qef{d7} in
$\Om\times(0,T)$. Then, $(\n,u)$ is said to be a global  entropy
weak solution if there exists some function $\Lambda\in
L^2(\Om\times(0,T))$ satisfying \qef{d3}, i.e.,
\bn\int_0^T\int \Lambda\vp dxdt
 &=&
 -\int_0^T\int\rho^{\al-1/2}\sqrt{\rho}u\vp_{x }dxdt\no&&
 - \frac{2\al}{2\al-1}
   \int_0^T\int(\rho^{\al-1/2})_x\sqrt{\rho}u\vp dxdt \lb{f50}
\en
 for any $\vp\in \Phi,$  and the  following uniform entropy
inequality holds\bma
 &\sup_{0\le t\le T}\int_\Om\left(|\sqrt{\n}u|^2
    +  |(\n^{\al-1/2})_x |^2
    + \pi(\n )\right)(x,t)\,dx\nnm\\
 &+\int_0^T\int_\Om
    ( |(\n^{(\ga+\al-1)/2})_x|^2
     +\Lambda^2)(x,t)\,dxdt
\nnm\\
\le
  & C_0\int_\Om (\mbox{$\frac{|m_0|^2}{\n_0}$}
    + |(\n_0^{\al-1/2})_x|^2
    + \pi_+(\n_0))(x)\,dx             \label{d4}
\ema with $C_0>0$ independent of  $T,$  and
 \bn
 \pi(\n)\triangleq\begin{cases}
    \n\log\n,&\mbox{  if }\quad \ga=1,\\
    \frac1{\ga-1}\,\n^\ga, &\mbox{  if }\quad \ga>1,
\end{cases}\quad
 \pi_+(\n)\triangleq\begin{cases}
    \n\log_+\n,&\mbox{  if }\quad \ga=1,\\
    \frac1{\ga-1}\,\n^\ga, &\mbox{  if }\quad \ga>1.
\end{cases}\lb{f51}
\en
\end{defi}

We have the following result on the existence of global entropy
weak solutions.\par

\opthm{(\textbf{Global existence})}\label{global_existence} Assume
that    \be
  \al> \frac{1}{2},\quad \ga> \frac{\al}{2}  .
   \label{l0e-ga}
\ee and that the initial data $(\n_0,m_0)$ satisfies \qef{d8} and
$\mbox{$\frac{|m_0|^{2+\nu}}{\n_0^{1+\nu}}$}\in L^1(\Om)$ for some
positive constant $\nu $.
Then for any $T>0,$ there exists a global   entropy weak solution
$(\n,u)$  to  the IBVP for the compressible Navier-Stokes
equations~\qef{l0a}-\qef{l0c} with boundary condition \qef{d6} or
\qef{d7} in $\Om\times (0,T)$ in the sense of
Definition~\ref{definition1}.

 Moreover, for the case of Dirichlet boundary condition \qef{d6}, if
$\al\in(1/2,3/2)$ and  $\nu$ satisfies \qef{c5} $($see Remark
\ref{nu-value} below$)$, then  in addition to \qef{s1}, the
solution $(\n,\sqrt{\rho}u)$ satisfies \be
  \sqrt{\rho}\,(\sqrt{\rho}{u})
  \in L^2(0,T;W_0^{1,(4+2\nu)/(4+\nu)}(\Om)),
 \lb{s2}
\ee
i.e., $\sqrt{\rho}\,(\sqrt{\rho}{u})$ satisfies the Dirichlet
boundary \qef{d6} in the sense of trace. \clthm

\begin{remark}\label{nu-value}
$(1)$. Theorem~\ref{global_existence} above holds for the
compressible shallow water equation  \qef{ii2}.

$(2)$. For the Dirichlet case \qef{d6},  one of the available ways
(available for the case $\al\in(1/2,3/2)$ and $\ga\ge 1$) to obtain
that $\sqrt{\n}\sqrt{\n}u\in L^p(0,T;W^{1,q} (\Om)),$  for some
$p\ge 1,q\ge 1,$ is to choose the positive constant $\nu$ in
Theorem~\ref{global_existence} such that \bn
 \begin{cases}
 \, \nu\in(0,2\ga-\al), &\mbox{for }
    \al\in (\mbox{$\frac{1}{2}$},1],\quad
   \ga\ge 1,
 \\[2mm]
 \, \nu\in\left[ \frac{2(\al+\ga-2)}{3-\al-\ga},
   \frac{2(2\ga-\al)}{1+\al-2\ga }\right],
 &\mbox{for }  \al\in(1,\mbox{$\frac{3}{2}$}),\quad
  \ga\in \left[1, \mbox{$\frac{1+\al}{2}$}\right),
  \\[2mm]
 \, \nu\in \left[\frac{4(\al-1)}{3-2\al},\infty\right),&\mbox{for }
    \al\in(1,\mbox{$\frac{3}{2}$}), \quad
  \ga\ge \mbox{$\frac{1+\al}{2}$}.
\end{cases}  \lb{c5}
\en
\end{remark}

Next, we show that there is a global  entropy weak solution
$(\n,u)$ in the sense of Definition~\ref{definition1} for which
the vacuum states and the structure of interface, if existing
initially, can be preserved for a short time, so long as the
initial data has additional regularity besides \ef{d8} and the
fluids and the vacuum states in initial data are connected
``smoothly". In addition, the weak solution $(\n,u)$ is actually a
unique regular solution for the short time. For simplicity, we
consider the case of one point vacuum state contained at
$x=x_0\in(0,1)$ in the initial data $(\n_0,m_0)=(\n_0,\n_0u_0)$
with additional regularity \bgr
    A_0|x-x_0|^{\sigma}
\le
  \n_0(x)
\le
  A_1|x-x_0|^{\sigma},
  \mbox{ for  any }  x \in \Om,
   \label{l0h}
\\
  u_0\in C^1(\bar\Om),\quad
  (\n_0^{\ga-1+1/2j})_x\in  L^{2j}(\Om),\
 \n_0^{-1+1/2j}(\n_0^{\al}u_{0x})_x
  \in L^{2j}(\Om),\   j=1,n,    \label{l0e}
\egr
with $n\ge 2 $ an  integer;  and in the case of continuous vacuum
state of one piece initially on $\Om^0=[x_0,x_1]\subset (0,1)$ in
the initial data, we require
 \bn
  \label{l0g}\begin{cases}
   A_0(x_0-x)^{\sigma} \le \n_0(x)
  \le A_1(x_0-x)^{\sigma}, & x\in[0,x_0),
\\
   \n_0(x) =0,
    m_0(x)=\n_0u_0(x)=0,& x\in[x_0,x_1],
\\
   B_0(x-x_1)^{\sigma}
    \le \n_0(x)\le B_1(x-x_1)^{\sigma},
   & x\in(x_1,1]\end{cases}\en and
\be\left\{\bln
 &(\n_0^{\ga-1+1/2j})_x\in  L^{2j}(\Om),\ j=1,n,\quad
 u_0\in C^1(\bar\Om\setminus\Om^0),
 \\[2mm]
 &\n_0^{-1+1/2j}(\n_0^{\al}u_{0x})_x
  \in    L^{2j}(\Om\setminus\Om^0),\ j=1,n, \label{l0ea}
 \eln\right.\ee
with  $n\ge 2 $ an integer. Here, $\sigma$, $A_0,\,A_1$, and
$B_0,\,B_1$ are positive constants, and the power
$\sigma\in(\sigma_-,\sigma_+)$ with positive constants
$\sigma_\pm$ given in \ef{beta} later. We also require that the
initial data $(\n_0,m_0)=(\n_0,\n_0u_0)$ given by \ef{l0c} is
consistent with boundary value for the Dirichlet boundary
condition.\par

We have the following results on short time structure of global
entropy weak solutions.
\opthm{(\textbf{Short time structure of vacuum states})}\label{short_time}
In addition to the assumptions of Theorem~\ref{global_existence},
assume further that
\be
  \al>\mbox{$\frac12$},\quad\gamma> \max\{1, \  \al\}
 \label{short-a}
 \ee
and that there is either one point vacuum state in initial data
$(\n_0,u_0)$ with \qef{l0h}--\qef{l0e} satisfied or a piece of
continuous vacuum states in initial data $(\n_0,u_0)$ with
\qef{l0g}--\qef{l0ea} satisfied.
Then, there exists a global  entropy weak solution $(\n,u)$ to the
IBVP for the compressible Navier-Stokes equations
\qef{l0a}-\qef{l0b} with initial data \qef{l0c} and boundary value
\qef{d6} or \qef{d7} in the sense of Definition \ref{definition1}.
\par

Moreover, there is a short time $T_*>0$, so that the global
entropy weak solution $(\n,u)$ is unique\footnote{Here the
uniqueness is specified for density $\n$ and momentum
$\n{u}=\sqrt{\n}\sqrt{\n}u$ for continuous vacuum states of one
piece.  } and regular on the domain $\Om\times[0,T_*]$, and the
initial structure of vacuum states is maintained for $t\in[0,T_*]$
in the following sense:
\par

For the case of one point vacuum state initially, \qef{l0h}, the
solution $(\n,u)$ is regular and unique on the domain
$\bar\Om\times[0,T_*],$ \bgr
 (\n,u)\in C^0(\bar\Om\times[0,T_*]),\quad
  u_x\in L^\infty(0,T_*;C^0(\bar\Om)),
                  \label{density}
\\
  \|u\|_{L^\infty(\bar\Om\times[0,T_*])}
 +\|u_x\|_{L^\infty([0,T_*];C^0(\bar\Om))}
 \le C(T_*). \label{velocity}
\egr
  The one point vacuum state
propagates along particle path, namely, there is one particle path
$x=X_0(t):[0,T_*]\rightarrow \ol\Om $ with $X_0(t)\in C([0,T_*] )$
defined by
 \be
 \dot{X}_0(t)=u(X_0(t),t),\quad X_0(0)=x_0\in (0,1), \label{point}
\ee
so that
\be
   a_-|x-X_0(t)|^{\sigma}
\le
  \n(x,t)
\le
  a_+|x-X_0(t)|^{\sigma}     \label{short-b}
\ee
for $(x,t)\in \Om\times[0,T_*]$, where the two positive constants
$a_\pm$ are independent of time $T_*$.
\par

In the case of a piece of continuous vacuum states initially,
\qef{l0g}, there are two particle pathes
$x=X_i(t):[0,T_*]\rightarrow \ol\Om $ with $X_i(t)\in C([0,T_*]
),i=0,1$ defined by
 \be
\dot{X}_i(t)=u(X_i(t),t),\quad X_i(0)=x_i\in (0,1),i=0,1,
\label{piece}
 \ee
so that it holds for some  positive constants $a_\pm,\,b_\pm$
independent of the time $T_*,$  \bgr
   a_-(X_0(t)-x)^{\sigma} \le \n(x,t)
\le
  a_+(X_0(t)-x)^{\sigma},     \label{short-c} \egr
   for  $(x,t)\in [0,X_0(t))\times[0,T_*],$ and \bgr
   b_-(x-X_1(t))^{\sigma}
    \le \n(x,t)\le b_+(x-X_1(t))^{\sigma}  \label{short-d}
\egr for  $(x,t)\in (X_1(t),1]\times[0,T_*]$ respectively, and the
interfaces separating  the fluid and vacuum   coincide with the
particle pathes \be
  \n(x,t) =0, \ \ \n{u}(x,t) =0,\quad
   (x,t)\in [X_0(t),X_1(t)]\times[0,T_*].    \label{short-e}
 \ee
  The solution $(\n,u)$ is regular and unique up to the vacuum
boundary
 \bgr
 \n\in C^0(\bar\Om\times[0,T_*]),\quad
  u\in C^0(\bar\Om\times[0,T_*]\setminus{\Om^0_{T_*}}),
                  \label{density-a}
\\
  \|u\|_{L^\infty(\bar\Om\times[0,T_*]\setminus{\Om^0_{T_*}})}
 +\|u_x\|_{L^\infty(\bar\Om\times[0,T_*]\setminus{\Om^0_{T_*}})}
 \le C(T_*). \label{velocity-a}
\egr
where $\Om^0_{T_*}=(X_0(t),X_1(t))\times[0,T_*]$.
\clthm

\oprem{} $(1).$ The constant exponents $\sigma_\pm$ are defined as
 $\sigma_\pm=\beta_\pm/(1-\beta_\pm)>0$ with $\beta_\pm$
determined by \be
 \beta_-
  =\max\{ \mbox{$\frac1{2\al}$},
   \ \mbox{$\frac1{\ga}$}(1-\mbox{$\frac{1}{2n}$})\},\quad
 \beta_+
  = \min\{1,\
    \mbox{$\frac1{\al}$}(1-\mbox{$\frac{1}{2n}$}),\
    \mbox{$\frac1{1+3\al}$}(4-\mbox{$\frac{1}{n}$})\},
    \label{beta}
\ee while the positive constants $a_\pm$ are independent of the
time $T_*$.
\par

$(2).$ The regularity assumptions \qef{l0h}--\qef{l0e} are
satisfied for the following initial data
\[
\n_0(x)=\mbox{$\frac12$}(A_0+A_1)(|x-x_0|^2)^{(\sigma_-+\sigma_+)/4},
 \quad u_0(x)=0,\quad x\in\Om,
\]
and the regularity assumptions \qef{l0g}--\qef{l0ea} are satisfied
for the initial data
\[
\n_0(x)=\begin{cases}
 \frac12(A_0+A_1)(x_0-x)^{(\sigma_-+\sigma_+)/2}, & x\in[0,x_0),
 \\
 0,   & x\in(x_0,x_1),\\
\frac12(B_0+B_1)(x-x_1)^{(\sigma_-+\sigma_+)/2}, & x\in(x_1,1],
\end{cases} \quad u_0(x)=0,\quad x\in\Om.
\]
\clrem
\par

Next, we prove that for any global  entropy weak solution $(\n,
u)$ to  the IBVP \qef{l0a}-\qef{l0c} together with boundary
condition \qef{d6} or \qef{d7} in the sense of
Definition~\ref{definition1}, even though in some cases that the
vacuum states may exist for some finite time, for instance, in the
cases as shown by Theorem~\ref{short_time}, any possible vacuum
state has to vanish within finite time after which the density is
always away from vacuum. Simultaneously, not only can the velocity
field be defined in terms of the density and momentum, and the
nonlinear diffusion is represented in terms of the density and
velocity, but also the global  entropy weak solution $(\n, u)$ is
shown to be a unique and strong solution after the vanishing of
vacuum states. We have the following result.\par

\opthm{(\textbf{Vanishing of vacuum  states})}\label{entropy_solution}
Assume that
\be
 \al>\frac12,\quad \ga\ge 1. \label{vanishing-vacuum}
\ee
 Let $(\n,u)$ be any global  entropy weak solution to  the
IBVP~\qef{l0a}-\qef{l0b} with initial data \qef{l0c} and boundary
value \qef{d6} or \qef{d7} in the sense of
Definition~\ref{definition1}. Then, there exist some time $T_0>0$
(depending on initial data) and a constant $\n_-$ so that \be
 \inf_{x\in\bar\Om}\n(x,t)\ge \n_->0,\quad t\ge T_0,  \label{vvs2}
 \ee
and the global entropy weak solution $(\n,u)$ becomes a unique
strong solution $(\n,u)$ for $t\ge T_0$ and satisfies
\bn \begin{cases}
 \n \in L^\infty(T_0,t;H^1(\Om)), \quad
 \n_t \in L^\infty(T_0,t;L^2(\Om)),
\\
 u\in H^1(T_0,t;L^2(\Om))\cap
      L^2(T_0,t;H^2(\Om)) ,  \end{cases} \label{vvs3b}
 \en
with velocity $u$ and   nonlinear diffusion term given by
 \be
{u}\triangleq\frac{\sqrt{\n}u}{\sqrt{\n}},\quad
 (\n^{\al } u_x)_x  =  \Lambda_x,
 \ee
respectively. In addition, for
 \[
  u_s\triangleq\begin{cases}0&\mbox{ for the Dirichlet case,}\\
\dis{\frac{1}{\bar{\n_0}}\int_\Om} m_0dx&\mbox{ for the periodic
case,}\end{cases}
 \]
there exist two positive constants  $\mu_0, c_0 $ both depending on
initial data $(\n_0,m_0)$ and $\n_-,$  such that\be
   \|(\n-\bar{\n_0},u-u_s)(\cdot,t)\|_{L^2( \Om)}
    \le
  c_0e^{-\mu_0(t-T_0)},\quad t>T_0,      \label{vvs4}
 \ee
where and what follows  $\ol f$  denotes the average of $f$ over the
bounded domain $\Om$, i.e.,
 \[
\overline{f}=\frac{1}{|\Omega|}\int_{\Omega}f(x)dx=\int_{\Omega}f(x)dx.
\]
\clthm

\begin{remark} $(1).$
 The Theorem \ref{entropy_solution} shows that any possible
vacuum states must vanish in finite time. This theory applies to
the compressible shallow water equation~\qef{ii2}.
\par
$(2).$ It is easy to verify (see the proof of Proposition \ref{tt2})
that the phenomena of vacuum vanishing \qef{vvs2} in finite time
actually happens for any global weak solution $(\n,u)$ to the
IBVP~\qef{l0a}-\qef{l0c} with boundary condition \qef{d6} or
\qef{d7} in the sense of Definition~\ref{definition} satisfying the
following entropy inequality \bma
 \lefteqn{\sup_{0\le t\le T}\int_\Om\left(|\sqrt{\n}u|^2
    +  |(\n^{\al-1/2})_x |^2
    + \pi(\n )\right)(x,t)\,dx +\int_0^T
      \|(\n^{(\ga+\al-1)/2})_x\|_{L^2}^2
        dt
}\no&& \le
  & C_0\int_\Om (\mbox{$\frac{|m_0|^2}{\n_0}$}
    + |(\n_0^{\al-1/2})_x|^2
    + \pi_+(\n_0))(x) dx      \aaa\aaa\aaa\aaa  \aaa\aaa
    \lb{d100}
\ema with $C_0$ independent of $T.$
\end{remark}

Finally, for any global  entropy weak solution $(\n,u)$ to the
IBVP~\ef{l0a}-\ef{l0c} together with boundary  condition \qef{d6}
or \qef{d7} in the sense of Definition~\ref{definition1}, the
density is continuous, i.e.,  $\n\in C(\ol\Om\times[0,T])$ for any
$ T>0$, due to \ef{s1} and \ef{d1}.  Thus, the continuity of $\n$
and Theorem~\ref{entropy_solution} imply that if the density
contains vacuum states at least at one point, then there exists
some critical time $T_1\in[0,T_0)$ with $T_0>0$ given by
\qef{vvs2} and a nonempty subset $\Om^0\subset\bar\Om$ such that
\bn\begin{cases}
 \n(x ,T_1)=0, &\forall \ x  \in\Om^0\\
 \n(x,T_1)>0,\ &\forall\ x\in\bar\Om\setminus\Om^0,\\
    \n(x,t)>0,  &
  \forall\  (x,t)\in\bar\Om\times(T_1,T_0].
  \end{cases}  \label{blowup-1c}
 \en
It follows from \ef{vvs3b} easily that for any $\delta>0,$ it
holds
 \be
 \int_{T_1+\de }^{T_0}\|u_x\|_{L^\infty}ds <\infty.
\ee
Under the condition that vacuum states appear, we shall  prove
that  the spatial derivative of velocity (if regular enough and
definable) blows up in finite time as the vacuum states vanish,
even if the solution is regular enough for short time so that the
velocity field and its derivatives are bounded as shown by
Theorem~\ref{short_time}.
\par
\opthm{(\textbf{Finite time blow-up})}\label{blow-up} Let $(\n,u)$
be any global   entropy weak solution, which contains vacuum
states at least at one point for some finite time,  to  the IBVP
for the compressible Navier-Stokes equations~\qef{l0a}-\qef{l0c}
with boundary condition \qef{d6} or \qef{d7} in the sense of
Definition~\ref{definition1}. Let $T_0>0$ and $T_1\in [0,T_0)$ be
the time such that \qef{vvs2} and \qef{blowup-1c} holds
respectively.

 Then, the solution $(\n,u)$ blows up as vacuum states vanish.
Namely, for $T_1$ satisfying \qef{blowup-1c} and for given  any
fixed $\eta>0,$ it holds
 \be
 \lim_{t\to T_{1}^+}\int^{T_1+\eta}_{t}\|u_x\|_{L^\infty}ds
=\infty.\label{blowup-1a}
\ee
On the other hand, if there exists some $T_2\in(0,T_0)$ such that
the weak solution $(\n,u)$ satisfies
\[
\|u\|_{ L^1(0,T_2;W^{1,\infty}(\Om))}<\infty,
\]
then, there is a time $T_3\in [T_2,T_0)$ so that the blowup
phenomena happens for  $(\n,u)$, i.e.,
 \be
 \lim_{t\to T_{3}^-}\int^{t}_{0}\|u_x\|_{L^\infty}ds
=\infty. \label{blowup-1g}
\ee \clthm

\oprem Theorem \ref{blow-up} implies that   for any global entropy
weak solution $(\n,u)$  to  the IBVP for the compressible
Navier-Stokes equations~\qef{l0a}-\qef{l0b} with initial data
\qef{l0c} and boundary value \qef{d6} or \qef{d7} in the sense of
Definition~\ref{definition1}, which contains vacuum states at least
at one point initially, the finite time blowup phenomena
\qef{blowup-1a} happens for such solution $(\n,u)$. \clrem

\oprem{}\label{whole_phenomena}
Theorems~\ref{global_existence}--\ref{blow-up} provide a complete
dynamical description on the vanishing of vacuum states and
blow-up phenomena  for the global entropy weak solutions to  the
compressible Navier-Stokes equations with density-dependent
viscosity. That is, a global  entropy weak solution exists for
general large initial data with finite entropy. For short time,
such weak solution is unique and regular with well-defined
velocity field subject to additional initial regularity, and any
existing vacuum state is maintained with the same interface
structure as initial. Then, within finite time the vacuum states
vanish definitely and the velocity blows up $($even if it is
regular enough and definable along the interfaces$)$. After the
vanishing of vacuum states, the global  entropy weak solution
becomes a strong one and tends to the non-vacuum equilibrium state
exponentially in time. This dynamical phenomena is quite similar
to those well-known for the 3-D incompressible Navier-Stokes
equations. However, before the time of vacuum-vanishing, the
uniqueness of the global   entropy weak solution to  the
compressible Navier-Stokes equations with density-dependent
viscosity subject to the initial data is not known yet. \clrem

\oprem{}\label{application1}
All theories established in
Theorems~\ref{global_existence}-\ref{blow-up} fit the shallow
water equation \qef{ii2}. We believe that such phenomena described
by Theorems~\ref{global_existence}-\ref{blow-up} are also observed
for other compressible fluids with density-dependent viscosity,
such as Navier-Stokes equations with capillarity and/or drag
friction, Navier-Stokes--Poisson system, etc. \clrem
\par

\oprem{}
It is interesting to investigate the (global) dynamics of
(one-dimensional) interface connecting vacuum from initial time
until the vanishing of vacuum in order to investigate the dynamics
of interface and the vanishing of vacuum state and to verify the
formation of singularity for general case.
\clrem

\oprem{}\label{extension} $(1)$. Another interesting problem is
whether the phenomena of the vanishing of vacuum states and blow-up
happens for multi-dimensional compressible isentropic Navier-Stokes
system with density-dependent viscosity, especially on spatial
bounded domain with the Dirichlet boundary condition~\qef{d6}.
 \par

$(2)$. It is also interesting to study whether any vacuum states
shall vanish within finite time and blow-up phenomena happens for
the multi-dimensional full Navier-Stokes system with
density-dependent viscosity. It is not obvious yet $($although we
expect$)$ since it is not clearly understood yet how the dynamics
of temperature and heat-conduction shall affect the global
existence of weak solutions and the evolution of vacuum states,
especially in the case of spherical symmetry under the Dirichlet
boundary condition \qef{d6}. This is under further
investigation~\cite{LiLiXin06b}. \clrem

\section{Global existence of entropy weak solutions}
\setcounter{equation}{0} \label{Global-Theory}

 In this section, we will establish the existence of global entropy
weak solutions to the IBVP for the compressible Navier-Stokes
equations \qef{l0a}-\qef{l0c} together with \qef{d6} or \qef{d7}.
Since the compactness arguments are straightforward in the framework
of Mellet-Vasseur~\cite{MelletVasseur05}, we only need to construct
a sequence of approximate solutions by using some ideas developed by
Jiang-Xin-Zhang in \cite{JiangXinZhang06} for one-dimensional
compressible Navier-Stokes equations with free boundary, to
establish some uniform a priori estimates, in particular, the lower
and upper bounds of the density for the approximate solution
sequence with the help of BD entropy, and to justify the boundary
condition \ef{s2} for (limiting) weak solution for the Dirichlet
case as follows.

\underline{\it Outline of proof of
Theorem~\ref{global_existence}}.\quad  {\it Step 1. Construction of
smooth approximate solutions}. Let us consider the following
approximate compressible Navier-Stokes equations inspired by
Jiang-Xin-Zhang in \cite{JiangXinZhang06}
 \bgr
 \rho _{\ve t}
+ ( \r    \uu )_ x  = 0, \label{l2a}
\\
 ( \r\uu  )_t
+ ( \r u_\ve^2 + p( \r ))_x
            -(\mu_\eps( \r )u _{\ve x})_x=0,
  \label{l2b}
  \\
 ( \r , \r\uu )(x,0)
=( \n_{ 0\ve}, m_{ 0\ve})(x),  \label{l2c}
 \egr
with one of the following boundary conditions
 \be
 \uu(0,t)=\uu(1,t)=0,\lb{a40}
 \ee
 or
 \be
 \r ,\uu  \mbox{ are periodic in }x \mbox{ of period }1.\lb{a41}
 \ee
The viscosity $\mu_\eps$ is given by
\be
 \mu_\eps(\n)=\n^\alpha + \eps \n^\theta,
  \quad \eps >0,\quad \theta \in (0,1/2),  \label{viscosity}
\ee
where we remark (see Remark~\ref{lower1} below) that although the
modified viscosities \ef{viscosity} are limited to the case $\theta
\in (0,1/2)$, the methods adapted here to construct approximate
solutions also can be applied to the case $\theta \in (0,1)$ and
$\gamma>\theta$ after some modification, the interesting reader can
refer to \cite{JiangXinZhang06,GuoJiuXin2007} and references
therein, we omit the details here. The initial data $ \n _{0\ve} , m
_{0\ve} \in C^\infty(\ol\Om)$ satisfies
  \be
\begin{cases}
\n _{0\ve}\to \n_0\mbox{ in }L^1(\Om),\
 \n _{0\ve} ^{\al-1/2}\to \n_0^{\al-1/2} \mbox{ in }H^1(\Om),\\
(m _{0\ve})^2(\n _{0\ve})^{ -1 }\to
  m_0 ^2 \rho _0 ^{ -1 },
|m _{0\ve}|^{2+\nu}(\n _{0\ve})^{ -1-\nu }\to
 | m_0| ^{2+\nu} \rho _0 ^{ -1-\nu }\mbox{ in }L^1(\Om)
\end{cases}\lb{e5}
 \ee
as $\ve\to 0^+$, and
  \be
  \n _{0\ve}\ge c_0\ve^{1/(2 \al-2\te)}    \lb{e4}
  \ee
for some $c_0$ independent of $\ve.$  By \ef{e5} we can assume
without the loss of generality that
 \be
\int_\Om \n_{0\ve}dx =1.\lb{a101}
 \ee
By the standard arguments (see \cite{JiangXinZhang06} for
reference), after applying the classical theory of parabolic and
hyperbolic equations, one can obtain that   there exists some
$T_*>0 $ such that the approximate problem \ef{l2a}-\ef{l2c}
together with the boundary condition \ef{a40} or \ef{a41} has a
unique smooth solution $(\r,\uu)$  on $[0,T_*]$ with the density
away from vacuum, i.e., \bn \r(x,t)>0,  \mbox{  for  all }
(x,t)\in
 \ol\Om\times[0,T_*].\lb{a100}
\en

\textit{Step 2. The a-priori estimates.} To extend the local
solution globally in time, one needs to control the lower and upper
bounds of the density and get some a-priori estimates. As mentioned
in the introduction, this relies on the BD entropy inequality
developed by Bresch and Desjardins (see
\cite{BreschDesjardinsLin05,BreschDesjardins03,
BreschDesjardins03b,BreschDesjardins05} for instance). In the
Eulerian coordinates, the BD entropy for the approximate solutions
$(\r,u_\eps)$ of the IBVP problem for the system \ef{l2a}--\ef{l2a}
reads (see \cite{BreschDesjardinsLin05,BreschDesjardins03,
BreschDesjardins03b,BreschDesjardins05} for instance) as
 \bn
  \lefteqn{\int_\Om(|\sqrt{\r}u_{\eps}|^2
    + |\mbox{$\frac{\mu(\r)_x}{\sqrt{\r}}$}|^2
    + |\mbox{$\frac{\eps((\r)^\theta)_x}{\sqrt{\r}}$}|^2
    + \pi(\r))(x,t)\,dx} \no
&&\aaa+\int_0^t\int_\Om (\mu_\ve({\r})|u_{\ve x}|^2
     +\mbox{$\frac{p'(\r)\mu_\ve'(\r)}{\r }$}
                |\rho_{\ve x}|^2)(x,s)\,dxds
 \le C,\   t>0
                \label{len-1dc3}
\en with positive constant $C$ independent of $t$ and $\ve.$ Indeed,
multiplying \ef{l2b} by $u_\eps$, integrating by parts on $[0,1]$,
and using \ef{l2a},  we get the classical entropy estimate
 \be
  \int_\Om(|\sqrt{\r}u_{\eps}|^2
    + \pi(\r))(x,t)\,dx
+\int_0^t\int_\Om \mu_\ve({\r})|u_{\ve x}|^2\,dxds \label{len-1dc3b}
 \le C,\   t>0
 \ee
with positive constant $C$ independent of $t$ and $\ve.$

Multiplying the transport equation \ef{l2a} by
$\varphi_\eps'(\rho_\eps)$ with $\n\varphi_\eps'(\n)=\mu_\ve(\n),$
i.e., $\varphi_\ve(\n )=\frac{1}{\al} \n ^\al+\frac\ve\theta
\n^\theta$ as in \cite{BreschDesjardinsMetivier07}, we get
\be
 \varphi_\eps(\rho_\eps)_t
 +  \varphi_\eps(\rho_\eps)_xu_\eps
 +\varphi_\eps'(\rho_\eps)\rho_\eps u_{\eps x}=0.
\ee
Differentiating above equation with respect to $x$ and denoting
$v_\eps=\varphi_\eps(\rho_\eps)_x/\rho_\eps$, we have
\be
 (\rho_\eps v_\eps)_t
 + (\rho_\eps v_\eps u_\eps)_x
 +(\mu_\eps (\rho_\eps) u_{\eps x})_x=0.\label{new-2}
\ee
Summing the equations \ef{new-2} and \ef{l2a} together, multiplying
by $u_\eps+v_\eps$ and integrating by parts on $[0,1]$, we have
after a straightforward computation that
 \bn
  \lefteqn{\int_\Om(|\sqrt{\r}u_{\eps}|^2
    + |\mbox{$\frac{\varphi_\eps(\rho_\eps)_x}{\sqrt{\r}}$}|^2
    + \pi(\r))(x,t)\,dx} \no
&&\aaa+\int_0^t\int_\Om (\mu_\ve({\r})|u_{\ve x}|^2
     +\mbox{$\frac{p'(\r)\varphi_\ve'(\r)}{\r }$}
                |\rho_{\ve x}|^2)(x,s)\,dxds
 \le C,\   t>0
 \en
with positive constant $C$ independent of $t$ and $\ve.$ This
together with \ef{len-1dc3b}
 \bma
(\r^{2\al-2}+\ve^2\r^{2\te-2})
 (\n_{\ve x})^2\le & (\varphi_\ve(\r)_x)^2
  = (\r^{2\al-2}+2\ve\r^{\al+\te-2}+\ve^2\r^{2\te-2})(\n_{\ve
  x})^2\nnm\\
 \le &  2(\r^{2\al-2}+\ve^2\r^{2\te-2})
 (\n_{\ve x})^2,\nnm
 \ema
gives rise to the BD entropy \ef{len-1dc3}. 

In order to derive the lower bound and upper bounds of density for
the approximate solution $(\r,\uu)$, it is convenient to make use of
the Lagrangian coordinates $(y,t) $ with
\[
y=\int_0^x \r (z,t)dz,t=t,
\]
because the approximate solution $(\r,\uu)$  is away from vacuum at
least for short time due to \ef{a100} and fluid density transports
along particle path. In Lagrangian coordinates, we have the
equivalent system for $(\r,\uu)$ as follows
 \bgr
 \n_ {\ve t }  +   \r  ^2 u_{\ve y}   = 0, \label{l5a}
\\
u_{\ve t}   + (p(\r ))_y
            -(\r \mu_\eps(\r ) u_{\ve y}  )_y=0,      \label{l5b}
\\
 ( \r , \uu)(y,0)
=( \n _{0\ve}, m_{0\ve} \n^{-1}_{0\ve})(y),  \label{l5c}
 \egr
together with the boundary conditions similar to \ef{a40} or
\ef{a41}, where $(y,t)\in  \Om_L\times(0,T_*]$ and
 \[
 \Om_L \triangleq  (0,L_{\ve} ), \mbox{ with } L_\ve
 =\int_\Om\n_{0\ve}(x)dx=1
  \]
due to \ef{a101}.

Correspondingly, the BD entropy equality \ef{len-1dc3} becomes in
Lagrangian coordinates as follows
 \bn
\lefteqn{
 \sup\limits_{0\le t\le T_*}E_\ve(t)
 + \int _0^{T_*}\int_\Ol\r \mu_\ve(\r ) (u_{\ve y} )^2dydt}\no&&
 +\int _0^{T_*}\int_\Ol
  \left(\r^{\ga+\al-2}+\ve\r^{\ga+\te -2}\right)(\n_{\ve y})^2 dydt
 \le
 C_0 E_\ve(0) \lb{e11}
 \en
with $C_0$ independent of $ {T_*}$ and $\ve>0$, and $E_\ve(t)$
defined as
 \[
 E_\ve(t)\triangleq\int_\Ol\left(u_{\ve}^2+((\r^\al)_{
y})^2+\ve^2((\r^\te)_{ y})^2+\r^{-1} {\pi}(\r )\right)dy.
 \]

Since
 \bn
 \int_\Ol\r^{-1}dy=\int_\Om
dx=|\Om|=1,\lb{e14}
 \en
it follows from the continuity of $\r$ that there exists some
$y_0(t)\in \Ol$ such that
 \[
 \r(y_0(t),t)=|\Ol|/|\Om|=1.
 \]
Hence, it implies from \ef{e11} that
 \bn
 \r^\al(y,t)&=&\r^\al(y_0(t),t)+\int^y_{y_0(t)}(\r^\al)_ydy\no
  &\le&
(|\Ol||\Om|^{-1})^{\al}+C_0E_\ve(0)+|\Ol|. \label{upper}
 \en
This yields the uniform upper bound (w.r.t. time $t\in[0,T_*]$)
for the density.

To obtain the lower bound for the density $\r $, we employ the idea
in \cite{JiangXinZhang06}. That is, we consider the upper bound for
$v_{\ve }= {\r^{-1} }, $  which can be estimated by \ef{e11} and
\ef{e14} as follow:
  \bma
   v_{ \eps}(y,t)
\le
 & \int_\Ol  v_{ \eps}(y,t) dy
   +\int_\Ol (v_{ \eps})^2| \n _{\ve y}|dy \nnm\\
 \le
 & 1
   +C\max_{y\in \ol\Om_L}(v_{ \eps})^{\theta+1/2}
    \|( \r  ^\theta)_y\|_{L^2(\Om)}
     (\int_\Ol v_{ \eps}(y,t) dy)^{1/2}\nnm\\
 \le
 & 1+\frac12\max_{y\in \ol\Om_L} v_{ \eps}
   +C(\ve^{-2}E_\eps(0))^{1/(1-2\te)}.   \label{lower}
\ema
 This shows that for all $(y,t)\in \ol\Om_L\times[0,T_*], $
 \bn
 \r (y,t)
\ge  \frac{\eps^{2/(1-2\theta)}}
            {2\eps^{2/(1-2\theta)}+C (E_\eps(0))^{1/(1-2\te)}}.
\label{len-1dc8}
  \en
Thus, we have obtained the uniform upper bound \ef{upper} and the
lower bound \ef{len-1dc8} (w.r.t. time $t\in[0,T_*]$) for density
$\r.$ Using these  uniform bounds, one can extend the local smooth
solution globally in time by standard arguments (see
\cite{JiangXinZhang06} for details). Moreover, the total mass of the
approximate solutions is conserved due to the boundary conditions
\ef{a40}
 \be
 \int_\Om \r (x,t)dx = \int_\Om\n_{0\ve} (x )dx,\quad
\forall\ t>0.\label{conservation}
 \ee
\begin{rem}\label{lower1}
Although the estimates \qef{lower} and \qef{len-1dc8} are
established for the approximate solutions in the case that the
modified viscosities \qef{viscosity} are limited to the case $\theta
\in (0,1/2)$, the methods adapted here can be applied to establish
the lower bounds like  \qef{len-1dc8} in the case $\theta \in (0,1)$
and $\gamma>\theta$. Indeed, one can also establish the following
estimates, similar to \cite{JiangXinZhang06,GuoJiuXin2007}, for the
approximate solutions $(\n,\uu)$ as
\[
 \sup_{t\in[0,T]}\int (\n^\theta)_x^{2n}(x,t)dx
 \le
 C(T)
\]
for any integer  $n\ge 1$ and $C(T)$ a constant independent of
$\eps,\theta$. Thus, we can verify similar to \qef{lower} that
\[\bln
   v_{ \eps}
\le
  1 + C\int  v^{1+\theta}_\eps|(\n^\theta)_x| dx
\le
   1+ C(T)(\int v_{ \eps}^{q(\theta+1)}dy)^{1/q}
    \le
  C + C(T)v^{(q(\theta+1)-1)/q}_\eps
\eln
\]
with $q=2n/(2n-1)$. This together with Young's inequality implies
the positive lower bound of density for integer $n$ sufficient large
but fixed so that $\theta<\frac{2n-1}{2n}$ so long as
$\theta\in(0,1)$.

\end{rem}

In addition to \ef{len-1dc3},  we show $\r^\eta\uu$ is actually
bounded in $L^\infty(0,T;L^{2+\nu}(\Om))$ for some $\eta\ge
1/(2+\nu)$
 with $\nu>0$  given by \bn
\begin{cases}\nu>0 \mbox{ arbitrary }&\mbox{ if }2\ga-\al\ge 1,\\
\nu\in(0,2(2\ga-\al)/(1+\al-2\ga)]&\mbox{ if }2\ga-\al\in (0,1).
 \end{cases}\lb{b4}\en
Note here that for $\al\in(1/2,3/2)$  the constant $\nu$ defined by
\qef{c5} actually satisfies \qef{b4}. It follows from \ef{len-1dc3}
that
 \be
 \sup\limits_{t\ge 0}
 \|(\r^{\al-1/2})_x \|_{L^2(\Om)}\le \ol C,
\ee which together with \ef{conservation} and  the mean value
theorem gives rise to
 \bn\label{a10}
 \sup\limits_{t\ge 0}\|\r\|_{L^\infty}^{\al-1/2}
\le
   \ol{\n_{0\ve}}^{\, \al-1/2}
 + \|(\r^{\al-1/2})_x \|_{L^2(\Om)}
\le
 \ol C
\en where $\ol C$ denotes some generic positive  constant depending
on $E_{\ve }(0)$ but independent of $T.$ With the help of  \ef{a10},
we can make use of the idea  due to \cite{MelletVasseur05} to derive
the following a-priori estimate  
 \bn
\sup\limits_{0\le t\le T}\int_\Om\r|\uu|^{2+\nu}dx\le C,\lb{a14}
 \en
which together with \ef{a10}
implies for $\eta(2+\nu)\ge 1$ that
 \bn
 \int_\Om\left(\r^\eta |\uu| \right)^{2+\nu }dx
 = 
 \int_\Om \r^{\eta(2+\nu)}|\uu|^{2+\nu}dx
 \le
 C\sup\limits_{0\le t\le T}\int_\Om\r|\uu|^{2+\nu}dx\le C,\lb{a15}
 \en
and, in particular, for $\eta=1/2$ that
 \be
  \sup\limits_{0\le t\le T}
  \int_\Om\left(\sqrt{\r }\,|\uu| \right)^{2+\nu }dx\le C  \lb{a15z}
 \ee
with $C$ one generic positive constant depending on both $T$ and
$E_{\ve }(0)$.

Thus, we can follow the compactness arguments in
\cite{MelletVasseur05} to prove the convergence of the approximate
solution $(\r,\uu)$ to some expected weak solution $(\rho,u)$ as
$\eps\to0+$. Namely, it holds
 \bgr
 \r\longrightarrow \rho \quad\mbox{ in }C([0,T]\times \ol\Om), \lb{a13}
 \\
  (\r^{\al-1/2})_x \rightharpoonup (\rho^{\al-1/2})_x
  \quad \mbox{ weakly in }L^2(\Om\times(0,T)),\lb{a29}
 \\
\sqrt{\r}\uu\longrightarrow \sqrt{\rho}u ,\quad
\r^\al\uu\longrightarrow \rho^{\al}u\quad\mbox{ in }
L^{2+\nu/2}(\Om\times(0,T)),
 \egr
and there exists some function $\Lambda\in L^2(\Om\times(0,T))$ such
that
 \bn
 \n^{\al }_\ve u_{\ve x}\rightharpoonup \Lambda
  \quad\mbox{ weakly in }  L^2(\Om\times(0,T)),
  \mbox{ as }\ve\to 0.                           \lb{f52}
 \en

{\it Step 3. Justification of the Dirichlet boundary condition
\ef{s2} for $\al\in(1/2,3/2) $.}  For $\al \in(1/2,1],$ it follows
from \ef{len-1dc3} and \ef{a15} that  \bn \lefteqn{ \| (\r
{u}_{\ve})_x \|_{L^2(0,T;L^{(4+2\nu)/(4+\nu)}(\Om))}}
 \no&&
\le
   C\left\|\r^{\al/2} u_{\ve x}\right \|_{L^2(0,T;L^2(\Om))}
 + C\left\|\left(\r^{\al-1/2}\right)_x\r^{1/2}\uu
     \right\|_{L^2(0,T;L^{(4+2\nu)/(4+\nu)}(\Om))}
 \no&&
\le
 C+C\left\|\left(\r^{\al-1/2}\right)_x\right\|_{L^\infty(0,T;L^2(\Om))}
\cdot\left\|\r^{1/2}\uu\right \|_{L^\infty(0,T;L^{2+\nu}(\Om))}
 \no&&
\le
 C.\lb{a16}
\en

For $\al\in(1,3/2),$  some additional estimates are needed here.
First, for   $\ga\ge (1+\al)/2,$ it is easy to check that for $\nu$
defined by \ef{c5} that $(3/2-\al)/(2+\nu)>1$.
Hence, for $\beta=(4+2\nu)/(4+\nu)$ one deduces from \ef{a15} and
\ef{a10} that \bn
 \lefteqn{
 \left\|\left(\r\uu\right)_x\right\|_{L^2(0,T;L^{\beta}(\Om))}}\no&&
\le
  C\left\|\r^{\al/2} u_{\ve x}\right\|_{L^2(0,T;L^2(\Om))}
 +C\left\|\left(\r^{\al-1/2}\right)_x
          \r^{ -\al+3/2}\uu\right\|_{L^2(0,T;L^{\beta}(\Om))} \no&&
\le
 C+C\left\|\left(\r^{\al-1/2}\right)_x\right\|_{L^\infty(0,T;L^2(\Om))}
 \left\|\r^{ -\al+3/2}\uu\right\|_{L^\infty(0,T;L^{2+\nu}(\Om))} \no&&
\le
 C.\lb{b7}
\en In the case that   $1\le \ga< (1+\al)/2,$ \ef{c5} yields that
$(3-\al-\ga)(2+\nu)\ge 2.$ Hence, one derives from \ef{a15} and
\ef{len-1dc3} that
 \bn
  \lefteqn{
  \left\|\left(\r\uu\right)_x\right\|_{L^2(0,T;L^{\beta}(\Om))}}\no&&
\le
  C\left\|\r^{\al/2} u_{\ve x}\right \|_{L^2(0,T;L^2(\Om))}
 +C\left\|\left(\r^{(\ga+\al-1)/2}\right)_x
           \r^{ (3-\al-\ga)/2}\uu\right\|_{L^2(0,T;L^{\beta}(\Om))} \no&&
\le
 C+C\left\|\left(\r^{(\ga+\al-1)/2}\right)_x \right\|_{L^2(0,T;L^2(\Om))}
\left\| \r^{( 3-\al-\ga
)/2}\uu\right\|_{L^\infty(0,T;L^{2+\nu}(\Om))}\no&& \le C. \lb{b8}
\en Thus, the estimates \ef{a40}, \ef{b7}   and \ef{b8} imply that
\ef{s2} still holds for   $\al\in(1,3/2).$ The proof of the Theorem
\ref{global_existence} on the existence of global entropy weak
solution  $(\rho,u)$ is completed.
\enddemo

\section{Dynamics of vacuum states for short time}
\setcounter{equation}{0}\label{Dynamics}
\subsection{Short time structure of vacuum states}
\label{short-time-structure}
We prove the Theorem~\ref{short_time} in this subsection in order
to study the short time structure of vacuum states for global
entropy weak solutions. To this end, it is sufficient to show that
there is a unique entropy  weak solution in short time for the
compressible Navier-Stokes equations~\ef{l0a}-\ef{l0b} with
initial data \ef{l0c} and boundary value \ef{d6} or \ef{d7} under
the assumptions of Theorem~\ref{short_time} as follows.
\par

\oppro{(\textbf{Vacuum states for short time})}\label{llocal-z}
Under the assumptions of Theorem~\ref{short_time}, there is a
short time $T'_*>0$ so that the unique entropy weak solution
$(\tilde{\n},\tilde{u})$ in the sense of
Definition~\ref{definition1} of the IBVP problem for the
compressible Navier-Stokes equations~\qef{l0a}--\qef{l0b} with
initial data \qef{l0c} and boundary condition \qef{d6} or \qef{d7}
exists on the domain $\Om\times[0,T'_*]$. The initial vacuum state
\qef{l0e} or \qef{l0ea} is also propagated for the short time
$t\in[0,T'_*]$, more precisely, the properties
\qef{density}--\qef{short-b} or \qef{piece}--\qef{velocity-a} hold
respectively for $(\tilde{\n},\tilde{u})$. In addition, it holds
that
 \be
 \|(\tilde{\n}^{\al-1/2})_x\|_{L^\infty(0,T'_*;L^2(\Om))}
 +
 \|((\tilde{\n}{\tilde{u})^2/\tilde{\n}},
   (\tilde{\n}\tilde{u})^{2+\nu}/\tilde{\n}^{1+\nu})
 \|_{L^\infty(0,T'_*;L^1(\Om))}
 \le
 C(T'_*). \label{short-g}
\ee
\clpro
 \demo
The proof of Proposition~\ref{llocal-z} will be completed  in
subsections~\ref{Lagrangian}--\ref{proof of proposition4.2} later.
\par
\bigskip
\underline{\it Proof of Theorem~\ref{short_time}}.\  This is a
consequence of  the Proposition~\ref{llocal-z} and
Theorem~\ref{global_existence}. In fact, the
Proposition~\ref{llocal-z} shows not only that there exists a
unique entropy weak solution $(\tilde{\n},\tilde{u})$ on the
domain $\Om\times (0,T'_*)$ in the sense of
Definition~\ref{definition1} to the IBVP problem for the
compressible Navier-Stokes equations~\qef{l0a}--\qef{l0b} with
initial data \qef{l0c} and boundary condition \qef{d6} or
\qef{d7}, but also that this short time entropy weak solution
satisfies all the properties \ef{density}--\ef{velocity-a}. Now,
choose time $T_*=T'_*-\delta$ with $\delta>0$ a constant small
enough. One can verify that at time $t=T_*$ the density and
momentum also satisfies the assumptions \ef{d8}, and particularly
\be
\begin{cases}
  \tilde{\n}(x,T_*)\ge 0 \mbox{ on }\Om,\quad
     \tilde{\n}\tilde{u}(.,T_*)=0,\
     \mbox{ on }\{x\in \Om\,|\,\n_0(x)=0\},
\\[2mm]
   \tilde{\n}(.,T_*)\in L^1(\Om),
   \quad (\tilde{\n}^{\al-1/2}(.,T_*))_x\in L^2(\Om),
\\[2mm]
  \mbox{$\frac{|\tilde{m}(.,T_*)|^2}{\tilde{\n}(x,T_*)}$}
  +\mbox{$\frac{|\tilde{m}(.,T_*)|^{2+\nu}}
               {\tilde{\n}^{1+\nu}(x,T_*)}$}  \in L^{1}(\Om).
\end{cases}                   \lb{d8z}
\ee
Thus, it follows from the Theorem~\ref{global_existence} that
there is a global entropy weak solution $(\hat{\n},\hat{u})$ for
time $t\ge T_*$ in the sense of Definition~\ref{definition1} to
the IBVP problem for the compressible Navier-Stokes
equations~\qef{l0a}--\qef{l0b} with initial data
\be
(\n,\n{u})(x,T_*)=(\tilde{\n},\tilde{\n}\tilde{u})(x,T_*),\quad
x\in \Om  \label{l0cz}
 \ee
and boundary condition \qef{d6} or \qef{d7}. Define $(\n,u)$  as
\be
(\n,u)=
\begin{cases}
 (\tilde{\n},\tilde{u}),
  \quad\mbox{ for \ } (x,t)\in\Om\times[0,T_*], \\
 (\hat{\n},\hat{u}),
  \quad \mbox{ for \ } (x,t)\in\Om\times[0,T_*].
\end{cases}  \label{global-a}
\ee It is easy to verify that $(\n,u)$ is a global entropy weak
solution in the sense of Definition~\ref{definition1} to the
compressible Navier-Stokes equations~\qef{l0a}--\qef{l0b} with
initial data \qef{l0c} and boundary condition \qef{d6} or \qef{d7}
under the assumptions of Theorem~\ref{short_time}. This global
entropy weak solution is actually unique and regular for short
time $t\in[0,T_*]$ and satisfies all the properties
\ef{density}--\ef{velocity-a}. The proof of
Theorem~\ref{short_time} is completed. \enddemo
\par

\subsection{Compressible Navier-Stokes in Lagrangian coordinates}
\label{Lagrangian}
In order to prove the Proposition~\ref{llocal-z},  we present an
equivalent proposition for the compressible Navier-Stokes equations
in the Lagrangian coordinates $(y,t)$, instead of the Eulerian
coordinates $(x,t)$, through the coordinate transformation \be
 y=\int_{0}^{x}\n(z,t)dz,
 \quad x\in(0,1), \ \ t\ge 0  \label{lag-dir}
\ee for the Dirichlet boundary condition, or \be
 y=\int_{X_0(t)}^{X_0(t)+x}\n(z,t)dz,
 \quad x\in(0,1), \ \ t\ge 0  \label{lag-dira}
\ee for the periodic boundary condition where $x=X_0(t)$ is a
particle path. In addition, both the case of initial one point
vacuum state and the case of a piece of initial continuous vacuum
states will be studied with some additional initial regularities.
\par

We first describe the equivalent proposition for the IBVP problem
for the compressible Navier-Stokes equations in the Lagrangian
coordinates in the case of initial one point vacuum state. Thus,
consider the compressible Navier-Stokes equations \bgr
 \n_t +  \n^2 u _y = 0, \label{l4a}
\\
 u_t  + p(\n)_y
            -(\n\mu(\n) u_y)_y=0,\quad y\in\Om,\ t>0
  \label{l4b}
 \egr
with initial data
\be
(\n,u)(y,0)=(\n_0(y),u_0(y)),\quad y\in \Om,   \label{l4c}
\ee
and the Dirichlet boundary condition
\be
 u(0,t)=u(1,t)=0,\quad t\ge 0  \label{l4bc-b}
\ee
or the periodic boundary condition
\be
  (\n,u)
  \mbox{\ is periodic w.r.t. $x$ of period one}. \label{l4bc-baa}
\ee
For the case of Dirichlet boundary, the initial data $(\n_0,u_0)$
is assumed to be consistent with the boundary values. Moreover, we
assume that there is one point vacuum state at $y=y_0$ for some
fixed point $y_0\in(0,1)$ and that the initial data is of
additional regularity
\be
\begin{cases}
 A_0|y-y_0|^{\beta} \le \n_0(y)
     \le A_1|y-y_0|^{\beta},\quad  y\in \Om,\\[2mm]
   u_0\in C^1(\bar\Om),\quad
  (\n^\ga_{0})_y,\
  (\n^{1+\al}_0u_{0y})_y\in L^{2}(\Om)\cap L^{2n}(\Om)
    \label{linit2}
\end{cases}
\ee with an integer $n\ge 2$, where $A_0,\,A_1$ are positive
constant, and the constant $\beta\in(\beta_-,\beta_+)$ with
$\beta_\pm$ determined by Remark~\ref{beta_pm}. Then the following
result holds in the case of initial point vacuum state.
\par
\oppro{(\textbf{Point vacuum state for short time})}\label{llocal}
Assume that \qef{short-a} and \qef{linit2}  hold.    Then, there
is a time $T'_*>0$ so that the unique regular weak solution
$(\n,u)$ with point vacuum of the IBVP \qef{l4a}--\qef{l4c} with
boundary condition \qef{l4bc-b} or \qef{l4bc-baa} exists on the
domain $\Om\times[0,T'_*]$ and satisfies \bgr
  \n\in C^0(\bar\Om\times[0,T'_*])\cap
           C^1([0,T'_*];L^2(\Om)),     \label{linit3a}
\\[2mm]
  u \in C^0(\bar\Om\times[0,T'_*])\cap
          C^1([0,T'_*];L^2(\Om)),      \label{linit3b}
\\[2mm]
  \n^{1+\al}u_y \in L^\infty(\Om\times[0,T'_*])
          \cap
          C^{1/2}([0,T'_*];L^2(\Om)),  \label{linit3c}
\\[2mm]
  \|(\n^\al)_y\|_{L^\infty([0,T'_*];L^2(\Om))}
 +\|\n{u_y}\|_{L^\infty(\Om\times[0,T'_*])}
 \le
  C(T_*).                                \label{linit3cz}
\egr
In addition, the initial point vacuum state is maintained for the
short time
\bgr
   \n(y_0,t)=0,\quad t\in[0,T'_*],
\\[2mm]
   a_-|y-y_0|^{\beta}\le \n(y,t)
\le   a_+|y-y_0|^{\beta},
  \  (y,t)\in[0,1]\times[0,T_*]. \label{l4bc-az}
\egr
Here, $a_\pm$ are positive constants independent of $T'_*$.
 \clpro
\demo The proof of the Proposition~\ref{llocal-z} will be given in
subsection~\ref{proof of proposition4.2}.
\par

\oprem{}\label{beta_pm}
The choice of $\beta_\pm>0$ is such that
\be
 \beta_- =\max\{ \mbox{$\frac1{2\al}$},
   \ \mbox{$\frac1{\ga}$}(1-\mbox{$\frac{1}{2n}$})\},\quad
 \beta_+ = \min\{1,\
    \mbox{$\frac1{\al}$}(1-\mbox{$\frac{1}{2n}$}),\
    \mbox{$\frac1{1+3\al}$}(4-\mbox{$\frac{1}{n}$})\}\label{beta-3}
\ee for integer $n\ge 2$. It should be emphasized here that all
the assumptions required here are satisfied for the shallow water
equations \ef{ii2}, i.e., $\ga=2,\ \al=1$.\par \clrem

Nest, we describe  the equivalent proposition for the IBVP problem
for the compressible Navier-Stokes equations in the Lagrangian
coordinates for the case of initial continuous vacuum of one
piece. In this case, we consider the IBVP problem for the
compressible Navier-Stokes equation~\ef{l4a}--\ef{l4c} with one of
following boundary conditions:
\begin{description}
\item{(1).} mixed boundary condition
 \be
\left\{\begin{aligned}
 & {u}(0,t)=\n(1,t)=0, \quad t\ge 0,\\
 &  A_0(1-y)^\beta \le \n_0(y)
     \le A_1(1-y)^\beta, \ y\in \Om, \label{l4bc-a}
\end{aligned}\right.
\ee
or \item{(2).} mixed free boundary condition
\be\left\{
   \begin{aligned}
   & \n(0,t)={u}(1,t)=0, \quad t\ge 0,\\
   & B_0y^\beta \le \n_0(y)
      \le   B_1y^\beta, \quad y\in\Om, \label{l4bc-c}
\end{aligned}\right.
\ee
 or \item{(3).}  free boundary condition
\be
\left\{
 \begin{aligned}
  & \n(0,t)=\n(1,t)=0, \  t\ge 0,\\
  & A'_0(y(1-y))^{\beta} \le \n_0(y)
     \le A'_1(y(1-y))^{\beta}, \ y\in \Om,\label{l4bc-d}
\end{aligned}\right.
\ee
\end{description}
where, $\beta\in(\beta_-,\beta_+)$, and $A_0,\,A_1$, $B_0,\,B_1$,
and $A'_0,\,A'_1$ are some given positive constants. The initial
data is also assumed to be regular \be
   u_0\in C^1(\bar\Om),\quad
  (\n^\ga_{0})_y,\
  (\n^{1+\al}_0u_{0y})_y\in L^{2}(\Om)\cap L^{2n}(\Om)\label{linit2z}
\ee with an integer $n\ge 2$. Then the following results for free
boundary value problems for the compressible Navier-Stokes
equations hold.
\oppro{(\textbf{Continuous vacuum state for short time})}
\label{llocal-zz}
Assume that \qef{short-a} and  \qef{linit2z} hold. Then, there is
a time $T'_*>0$ so that the unique regular weak solution $(\n,u)$
of the IBVP problem for the compressible Navier-Stokes
equations~\qef{l4a}--\qef{l4b} with initial data \qef{l4c} and
free boundary condition \qef{l4bc-a}, or \qef{l4bc-c}, or
\qef{l4bc-d} exists on the domain $\Om\times[0,T'_*]$ and
satisfies
\bgr
  \n\in C^0(\bar\Om\times[0,T'_*])\cap
           C^1([0,T'_*];L^2(\Om)),   \nnm
\\[2mm]
  u \in C^0(\bar\Om\times[0,T'_*])\cap
          C^1([0,T'_*];L^2(\Om)),   \nnm%
\\[2mm]
  \n^{1+\al}u_y \in L^\infty(\Om\times[0,T'_*])
          \cap
          C^{1/2}([0,T'_*];L^2(\Om)), \nnm
\egr
and
\be
  \|(\n^\al)_y\|_{L^\infty([0,T'_*];L^2(\Om))}
 +\|\n{u_y}\|_{L^\infty(\Om\times[0,T'_*])}
 \le
  C(T'_*).                            \nnm
\ee In addition, the initial vacuum state is also maintained for
short time, namely, it holds that \be
   a_-(1-y)^\beta \le \n(y,t)
\le   a_+(1-y)^\beta,
 \quad (y,t)\in(0,1)\times[0,T'_*] \nnm
\ee
corresponding to the mixed free boundary conditions  \qef{l4bc-a},
or
\be
   b_-y^\beta \le \n(y,t)
\le   b_+y^\beta,
 \quad (y,t)\in(0,1)\times[0,T'_*] \nnm
\ee
corresponding to the mixed free boundary conditions \qef{l4bc-c},
or
\be
   c_-(y(1-y))^{\beta}
   \le
\n(y,t)
   \le
   c_+(y(1-y))^{\beta},
  \quad (y,t)\in (0,1)\times[0,T'_*]
   \nnm
\ee
corresponding to the free boundary \qef{l4bc-d}. Here, $a_\pm$,
$b_\pm$ and $c_\pm$ are positive constants independent of $T'_*$.
\clpro
 \par
\demo
The short time existence, uniqueness and regularity of weak
solutions of the free boundary problems for the compressible
Navier-Stokes equations are well-investigated by many authors
(see~\cite{LiuXinYang98,LuoXinYang00,
JiangXinZhang06,Yangzhao02,Yangzhu02}).
Proposition~\ref{llocal-zz} can be proved in a similar way
as~\cite{LiuXinYang98,LuoXinYang00, JiangXinZhang06,Yangzhao02},
so we omit the details.
\par
\oprem{} Based on Proposition~\ref{llocal-zz} for the compressible
Navier-Stokes equations \qef{l4a}--\qef{l4c} with either mixed
free boundary condition \qef{l4bc-a} and \qef{l4bc-c} or free
boundary  \qef{l4bc-d} and the coordinates transformation from the
Lagrangian coordinates to the Eulerian coordinates, one can study
easily the IBVP problem for compressible Navier-Stokes equations
\qef{l0a}--\qef{l0b} in the Eulerian coordinates with either the
Dirichlet boundary condition \ef{d6} or  the periodic boundary
condition \qef{d7} in the case of initial continuous vacuum state
of one piece \qef{l0g}.

In fact, for the compressible Navier-Stokes equations
\qef{l0a}--\qef{l0b} with Dirichlet boundary condition \qef{d6} in
the Eulerian coordinates,  the short time existence of unique
solution subject to the case of a piece of initial continuous
vacuum state \qef{l0g}--\qef{l0ea} in initial data can be
constructed, in the Lagrangian coordinates, by combining two mixed
free boundary value problems \qef{l4bc-a} and \qef{l4bc-c}
together with one continuous vacuum state in-between as follows.
Denote two particle pathes $x=X_i(t)$ (assumed to be definable for
short time) starting from the initial vacuum boundary $x=x_0$ and
$x=x_1$ respectively as \be
 \dot{X}_i(t)=u(X_i(t),t),\quad X_i(0)=x_i,\ i=0,1,\nnm
\ee along which the vacuum boundary moves in the Eulerian
coordinates so that \be\begin{cases}
 (\n,\n{u})(x,t)=0,\quad &x\in[X_0(t),X_1(t)],\  t\ge 0,\\[2mm]
 \n{u}(x,t)>0,\quad &x\in[0,X_0(t))\cup(X_1(t),1],\  t\ge 0.
   \label{trans-f0}
\end{cases}
\ee

We first choose the coordinate transformation from the Eulerian
coordinates to the Lagrangian coordinates as \be
 \begin{cases}
 y=\int_0^{x}\n(z,t)dz,\quad
  x\in[0, X_0(t)],
\\[3mm]
 y_0=\int_0^{X_0(t)}\n(z,t)dz
    =\int_0^{x_0}\n_0(z)dz<1, \quad \mbox{conservation of mass},
          \label{trans-f1}
\end{cases}
\ee which gives $y\in[0,y_0]$ and the mixed free boundary
conditions \qef{l4bc-a}. The application of the
Proposition~\ref{llocal-zz} with $\Om$ replaced by $(0,y_0)$
implies, via the inverse coordinate transformation $x
=\int_0^{y}\n^{-1}(z,t)dz$, $y\in[0,y_0]$, the existence of unique
solution $(\n_l,u_l)$ of the compressible Navier-Stokes
equations~\qef{l0a}--\qef{l0b} on $[0,X_0(t)]\times[0,T_*]$ in the
Eulerian coordinates with the initial data
\[
(\n_l,u_l)(x,0)=(\n_0,u_0)(x),\quad x\in(0,x_0),\quad \n_0(x_0)=0,
\]
and mixed free boundary conditions
\[
u(0,t)=0,\quad \n(X_0(t),t)=0,\quad t\in[0,T_*].
\]

Next, we choose the coordinate transformation from the Eulerian
coordinates to the Lagrangian coordinates as \be
 \begin{cases}
 y= 1-\int_{x}^1\n(z,t)dz,\quad
  x\in[X_1(t),1],
\\[3mm]
 y_1=\int^1_{X_1(t)}\n(z,t)dz
    =\int^1_{x_1}\n_0(z)dz<1   \quad \mbox{conservation of mass},
            \label{trans-f2}
\end{cases}
\ee which gives $y\in[1-y_1,1]$ and the mixed free boundary
conditions  \qef{l4bc-c}. Then Proposition~\ref{llocal-zz}
implies, via the inverse coordinate transformation $x
=1-\int_y^{1}\n^{-1}(z,t)dz$, $y\in[1-y_1,1]$, the existence of
unique $(\n_r,u_r)$ of the compressible Navier-Stokes
equations~\qef{l0a}--\qef{l0b} on $[X_1(t),1]\times[0,T_*]$ in the
Eulerian coordinates
 with initial
data
\[
(\n_r,u_r)(x,0)=(\n_0,u_0)(x),\quad x\in(x_1,1),\quad \n_0(x_1)=0,
\]
and mixed free boundary conditions
\[
u(1,t)=0,\quad \n(X_1(t),t)=0,\quad t\in[0,T_*].
\]
Consequently , we can construct the short time unique solution
$(\n,u)$ to the IBVP problem for the compressible Navier-Stokes
equations~\qef{l0a}--\qef{l0c} with the Dirichlet boundary
condition \qef{d6} and a piece of continuous vacuum
state~\qef{l0g} in the initial data as
\[
(\n,\n{u})=\begin{cases}
 (\n_l,\n_l u_l), & \mbox{on}\quad [0,X_0(t)]\times[0,T_*],\\
 (0,0), & \mbox{on}\quad  (X_0(t),X_1(t))\times[0,T_*],\\
 (\n_r,\n_ru_r), & \mbox{on}\quad [X_1(t),1]\times[0,T_*].\\
  \end{cases}
\]

Similarly,  we can obtain the short time existence of unique
solution to the IBVP problem for the compressible Navier-Stokes
equations~\qef{l0a}--\qef{l0b} with periodic boundary condition
\qef{d7} which can be viewed as a free boundary problem after the
choice of the spatial reference point. The details will be
omitted. \clrem

\subsection{Proof of Propositions~\ref{llocal-z}--\ref{llocal}}
\label{proof of proposition4.2}
We first prove the Proposition~\ref{llocal} in this subsection
with the help of the a-priori estimates for (regularized)
solutions and the construction of approximate solutions by a
finite difference scheme, due to the modification of the ideas
used in
\cite{JiangXinZhang06,LuoXinYang00,OkadaMatsusuMakino02,Yangzhao02}.
Without the loss of generality, we only prove
Proposition~\ref{llocal} in the case of the Dirichlet boundary
conditions with one point vacuum state in the initial data.

First, we can easily derive some identities for (regularized)
solutions as in \cite{JiangXinZhang06,Yangzhao02}.

\oplem{}\label{lemma2}
Let $T>0$ and assume that the solution $(\n,u)$  of the IBVP
problem~\qef{l4a}--\qef{l4bc-b} exists for $t\in[0,T]$ with
$\n(y_0,t)=0$. Then, under the assumptions of
Proposition~\ref{llocal}, it holds that \bgr
\n^{1+\alpha}u_y(y,t)= \int_{y_0}^y u_t(z,t)dz  + \n^\gamma(y,t)
 = -\int_y^{y_0} u_t(z,t)dz  + \n^\gamma(y,t),        \label{leq2}
\\
\bln
   \n^{\alpha}(y,t) + \alpha\int_0^t\n^\gamma(y,s)ds
 = &\n_{0}^{\alpha}(y)
    +\alpha\int_0^t\int_y^{y_0} u_t(z,s)dzds \\
 = &\n_0^{\alpha}(y)
    -\alpha\int_0^t\int_{y_0}^y u_t(z,s)dzds.\label{leq3}
\eln
\egr
for all $y\in\Om$ with $y\neq y_0$.
\cllem

We also have the following useful a-priori estimates, whose proofs
are similar to those for  \ef{len-1dc3} and \ef{e11} and thus will
be omitted. \oplem{}\label{lemma3}
Let $T>0$ and assume that the solution $(\n,u)$  of the IBVP
problem~\qef{l4a}--\qef{l4bc-b} exists for $t\in[0,T]$. Then,
under the assumptions of Proposition~\ref{llocal}, it holds that
\bgr \bln
 & \|u(t)\|_{L^2(\Om)}^2+\|\n(t)\|^{\gamma-1}_{L^{\gamma-1}(\Om)}
  +\int_0^t\int |(\n^{(\ga+\al)/2})_y|^2dxds
 \le
   C_0,           \nnm
  \quad
   t\in[0,T]
\eln\\
    \|\n(t)\|_{L^\infty({\Om})}
   +\|(\n^\al(t))_y\|_{L^2(\Om)}
 \le C_0,\quad
  t\in[0,T] \label{leq3z}
\egr
with $C_0>0$ a constant.
\cllem

\oplem{}\label{lemma4}
Let $T>0$  and assume that the solution $(\n,u)$  of the IBVP
problem~\qef{l4a}--\qef{l4bc-b} exists for $t\in[0,T]$ with
$\n(y_0,t)=0$.  Then, under the assumptions of
Proposition~\ref{llocal}, there is a time $T'_*\in (0,T]$
$($depending on initial data$)$ so that it holds that
\be
   \int u_t^{2j}(y,t)dy
  +\int_0^t\int  u_t^{2j-2} \n^{1+\alpha}u_{yt}^2dxds
  \le
  C_1         \label{les1}
\ee
uniformly for $t\in[0,T'_*]$ with $j= 1,n$ and $C_1>0$ a positive
constant, and  that $\n{u}_y(t)\in C^0(\bar\Om)$ is uniformly
bounded for any $t\in[0,T'_*]$
\be
 \|\n{u}_y\|_{L^\infty(\Om\times[0,T'_*])}\le C_2    \label{les2}
\ee
with $C_2>0$ a constant. Moreover, the solution $(\n,u)$ is
continuous
 \be
    (\n,u)\in C^0(\bar\Om\times[0,T'_*]),
    \label{les2z}
\ee
and the initial one point vacuum state is maintained
\be
 a_-|y-y_0|^{\beta}\le \n(y,t)\le a_+|y-y_0|^{\beta},
 \quad (y,t)\in\Om\times[0,T'_*]     \label{vacuum-point}
\ee
with $a_+> a_->0$ two constants independent of time $T_*$.
\cllem
\demo
Let us first assume that the weak solution is regular enough so
that we can differentiate it through the equations and the
interface as in \cite{Yangzhao02}, and the density is of the form
\be
 a_*|y-y_0|^{\beta}
 \le \n(y,t)
 \le a^*|y-y_0|^{\beta},
   \quad (y,t)\in\Om\times[0,T]   \label{lasp2}
\ee with  $a_*,\,a^*$ two positive constants to be determined
later. It will be assumed further that the following a-priori
estimate holds \be
 \|\n{u}_y\|_{L^\infty({\Om\times[0,T]})}\le M_0\label{lasp1}
\ee for some positive $M_0$ to be determined later. We will prove
\ef{les1} only for $j=n$ below since the other case can be treated
similarly. Taking inner product between $\ef{l4b}_t$ and
$2n(u_t)^{2n-1}$ over $\Om$ leads to \be\bln
   \frac{d}{dt}\int |u_t|^{2n}dy
 + 2n\int (\n^\gamma)_{yt}(u_t)^{2n-1}dy
 - 2n\int (\n^{1+\alpha}u_y)_{yt} (u_t)^{2n-1}dy
 =0\nnm
\eln \ee from which, one deduces after integration by parts and
using Eq.~\ef{l4a} that \bma
   &\frac{d}{dt}\int |u_t|^{2n}dy
  + 2n(2n-1)\int  \n^{1+\alpha}u_{yt}^2(u_t)^{2n-2} dy \nnm\\
 =
 & 2n(2n-1)\gamma\int \n^{\gamma-1}\n_{t}u_{yt}(u_t)^{2n-1}dy\nnm\\
 &+ 2n(2n-1)(1+\alpha)\int \n^{\alpha}\n_tu_y u_{yt}(u_t)^{2n-2}dy\nnm\\
\le
 &  n(2n-1) \int  \n^{1+\alpha}u_{yt}^2(u_t)^{2n-2} dy\nnm
   +CC_0^{2\gamma-1-\alpha}M_0^2(1+\|u_t(t)\|_{L^{2n-2}}^{2n+2})\\
 & +C\int \n^{3+\alpha}u_y^4(u_t)^{2n-2}dy\label{lin1e}
\ema where one has used \ef{short-a}, \ef{leq3z}, the a-priori
assumption \ef{lasp1}, and H\"older's inequality. The last term on
the right hand side of \ef{lin1e} can be estimated as follows. For
the case $ \alpha\ge 1$, it follows from \ef{leq3z} and \ef{lasp1}
that
 \bma \int \n^{3+\alpha}u_y^4(u_t)^{2n-2}dy
 \le
  &\|\n^{3+\alpha}u_y^4\|_{L^\infty}
    \|u_t(t)\|_{L^{2n-2}}^{2n-2}
 \le C_0^{\alpha-1}M_0^4\|u_t(t)\|_{L^{2n-2}}^{2n-2}\nnm\\
 \le
  & CM_0^4(1+\|u_t(t)\|_{L^{2n}}^{2n+2}). \label{lin1zz}
\ema
For the case $\alpha\in (\frac12,1)$, since it holds by
\ef{leq2} that
\bma
  |\n^{1+\alpha}u_y | \le |\int_{y_0}^y u_t(z,t)dz|  + \n^\gamma
\le
 \|u_t\|_{L^{2n}}\cdot |y-y_0|^{(2n-1)/2n} +\n^\ga,
\ema
which together \ef{lasp2} implies
\bma
  |\n^{(3+\alpha)n}u_y^{4n}|
 =
 & |\n^{-(1+3\alpha)n}(\n^{1+\al}u_y)^{4n}|
 =
  \n^{-(1+3\alpha)n}
   \left[\int_{y_0}^y u_t(z,t)dz + \n^\gamma\right]^{4n}\nnm\\
\le
 & C \n^{-(1+3\alpha)n}
    (\|u_t\|_{L^{2n}}^{4n}\cdot |y-y_0|^{2(2n-1)}
     +\n^{4n\ga})\nnm\\
\le
 & C\|u_t\|_{L^{2n}}^{4n}\cdot
     |y-y_0|^{-(1+3\alpha)n\beta +2(2n-1)}
     +\n^{3n(\ga-\al) + n(\ga-1)},\label{linh}
\ema
the last term in \ef{lin1e} is estimated by
\bma
 & |\int\n^{3+\alpha}u_y^4(u_t)^{2n-2}dy|
\le
  \|u_t\|_{L^{2n}}^{2n-2}\cdot
  \left(\int\n^{(3+\alpha)n}u_y^{4n}dy\right)^{1/n}\nnm\\
\le
 &C \|u_t\|_{L^{2n}}^{2n-2}\cdot
  \left(\|u_t\|_{L^{2n}}^{4n}
       \int |y-y_0|^{-(1+3\alpha)n\beta +2(2n-1)}dy
     +\int \n^{3n(\ga-\al) + n(\ga-1)}dy\right)^{1/n}\nnm\\
\le
 &C (\|u_t\|_{L^{2n}}^{2n+2} +1).  \label{linhz}
\ema where one has used the fact $\beta <\frac{4n-1}{n(1+3\al)}$
due to \ef{beta-3}. Substituting \ef{lin1zz} and \ef{linhz} into
\ef{lin1e} shows \bma
   &\frac{d}{dt}\int |u_t|^{2n}dy
  + n(2n-1)\int  \n^{1+\alpha}u_{yt}^2(u_t)^{2n-2} dy \nnm\\
\le
 & CM_0^2(1+M_0^2)(1+\|u_t(t)\|_{L^{2n-2}}^{2n+2})
\le
  C(1+M_0^4)
     (\|u_t(t)\|_{L^{2n}}^{2n+2} +1).   \label{lin1z}
\ema

Set
 \be
 T_a=\min \{\mbox{$ \frac{2^n-1}{n2^nC(1+M_0^4)C_1},\quad
 \frac{\|u_t(0)\|^{2n}_{L^{2n}}}{C(1+M_0^4)}$},\quad T
  \}.\label{T_a}
\ee One can apply the Gr\"onwall's Lemma to obtain \ef{les1} for
$t\in [0,T_a]$ with $C_1$ given by \be
    C_1
  =: 2\|u_t(0)\|^{2n}_{L^{2n}}
  \ge
   \|u_t(0)\|^{2n}_{L^{2n}}+CT_a(1+M_0^4).           \label{C_1}
\ee
\par

To prove \ef{les2} and ensure the a-priori assumption~\ef{lasp1},
we use the equality \ef{leq2} to get that near $y=y_0$, it holds
that \be
  \n{u}_y(y,t)
 =   \n^{-\alpha}(y,t)\int_{y_0}^y u_t(z,t)dz
   + \n^{\gamma-\alpha}(y,t),\quad y\neq y_0,       \label{leq2a}
\ee which, together with \ef{short-a}, \ef{leq3z}, and the fact
$\beta \le \frac1\al\left(1- \frac{1}{2n}\right)$ due to
\ef{beta-3}, implies \be\bln
  \|\n{u}_y\|_{L^\infty_{\Om\times[0,t]}}
\le
 &  C_0^{\gamma-\alpha}
  + |\n^{-\alpha}(y,t)\int_{y_0}^y u_t(z,t)dz
    |_{L^\infty_{\Om\times[0,t]}}   \\
\le
 &   C_0^{\gamma-\alpha}
  + a_*^{-\alpha}|y-y_0|^{1-1/2n-\al\beta}\|u_t(t)\|_{L^{2n}}\\
\le
 & C_0^{\gamma-\alpha}+  a_*^{-\alpha}C_1^{1/2n}
 \le
   C_0^{\gamma-\alpha}+  2a_*^{-\alpha}\|u_t(0)\|_{L^{2n}}\\
 =
 &:C_2\le M_0              \nnm
\eln
\ee
so long as the constant $M_0$ is chosen as
\be
 M_0= 1 + C_0^{\gamma-\alpha}
        +  2a_*^{-\alpha}\|u_t(0)\|_{L^{2n}}.  \label{M_0}
\ee

Next, we verify the a-priori assumption~\ef{lasp2}. Set \be
  T_b=\min\{\,\mbox{$\frac{A_0^{\alpha}}{3{\al}C_1^{1/2n}}$},
    \quad T_a\}. \label{T_b}
\ee
It follows from the equation \ef{leq3} for $y\neq y_0$ and
$t\in[0,T_b]$ that
\be
\bln
  & \n^{\alpha}(y,t) + \alpha\int_0^t\n^\gamma(y,s)ds
 =
 \n_0^{\alpha}(y) +\alpha\int_0^t\int_{y_0}^y u_t(z,s)dzds\\
\ge
 & \frac23A_0^{\alpha}|y-y_0|^{\al\beta}
  + |y-y_0|^{\al\beta}
    (\frac13A_0^{\alpha}- C_1^{1/2n}\al T_b|y-y_0|^{1-1/2n-\al\beta})\\
\ge
 & \frac23A_0^{\alpha}|y-y_0|^{\al\beta} \eln           \label{leq3a}
\ee where we have used the facts $|y-y_0|\le 1$ and $\beta \le
\frac1\al\left(1- \frac{1}{2n}\right)$ due to \ef{beta-3}. On the
other hand, it follows from \ef{leq3} that \be \bln
  & \n^{\alpha}(y,t) + \alpha\int_0^t\n^\gamma(y,s)ds
 =
 \n_0^{\alpha}(y) +\alpha\int_0^t\int_{y_0}^y u_t(z,s)dzds\\
\le
 & A_1^{\alpha}|y-y_0|^{\al\beta}
  + \al t|y-y_0|^{1-1/2n}\|u_t(t)\|_{L^{2n}}\\
\le
 & (A_1^{\alpha}+ C_1^{1/2n}\al t)|y-y_0|^{\al\beta}.
\eln                 \label{leq3c} \ee
Define $T_c\in(0,T_b]$ by
\be
 T_c
  =\min\{\, \frac{A_0^{\alpha}}
            {3\al C(A_1^{\alpha}+ C_1^{1/2n}\al T_b)^{\ga/\al}},
         \quad T_a,\, T_b\}. \label{T_c}
\ee
Set

\[
Z(t)= \int_0^t\n^\gamma(y,s)ds.
\]
It follows from  \ef{leq3c} that \be
 (Z'(t))^{\al/\ga}
  \le
   (A_1^{\alpha}+ C_1^{1/2n}\al t)|y-y_0|^{\al\beta}\nnm
\ee
which implies for $t\in[0,T_c]$ that
\be
\bln
  \int_0^t\n^\gamma(y,s)ds
\le &
 T_c (A_1^{\alpha}+ C_1^{1/2n}\al T_c)^{\ga/\al}
  |y-y_0|^{(\ga-\al)\beta}|y-y_0|^{\al\beta}\\
\le &
 CT_c(A_1^{\alpha}+ C_1^{1/2n}\al T_b)^{\ga/\al}|y-y_0|^{\al\beta}.
    \label{leq3d}
\eln \ee As a consequence of \ef{leq3a}, \ef{leq3d}, and \ef{T_c},
one gets \be
 \n^{\alpha}(y,t)
 \ge\,
 \frac23A_0^{\alpha}|y-y_0|^{\al\beta} - \al\int_0^t\n^\gamma(y,s)ds
 \ge\,
 \frac13A_0^{\alpha}|y-y_0|^{\al\beta}.    \label{leq3f}
\ee From \ef{leq3c} and \ef{leq3f}, we can verify the a-priori
assumption~\qef{lasp2} and justify the property \ef{vacuum-point}
by simply choosing \be
 a_*=a_-= (A_0^{\alpha}/3)^{1/\al},\quad
 a^*=a_+=(A_1^{\alpha}+ C_1^{1/2n}\al T_c)^{1/\al} \label{leq3s}
\ee for any $t\in[0,T'_*]$ with time $T'_*=T_c$ determined by
\ef{T_c}. One then derives from Eq.~\ef{l4a}, \ef{leq3z}, and
\ef{les2} that \be
 \n^\alpha\in L^\infty(0,T'_*,H^1(\Om)),\quad
 (\n^\alpha)_t\in L^\infty(0,T'_*,L^2(\Om)).  \label{leq3za}
\ee while \ef{les2}, \ef{vacuum-point}, boundary condition
\ef{l4bc-b}, and \ef{les1} for $j=1$ imply that
 \be
u\in L^\infty(0,T'_*,W_0^{1,p}(\Om)),\quad
 u_t\in L^\infty(0,T'_*,L^2(\Om))              \label{leq3zb}
\ee for any $p\in(1,\beta^{-1})$, where one has used the fact
$\beta^{-1}>\beta_+^{-1}\ge1$ so that \be
\sup_{t\in[0,T_*]}\|u_y\|_{L^p(\Om)}
 =\sup_{t\in[0,T_*]}\|\n^{-p}\|_{L^1(\Om)}
  \cdot\|\n{u}_y\|^p_{L^\infty(\Om\times[0,T'_*])}  \label{leq3zc}
 \le C.
\ee \ef{leq3za}--\ef{leq3zc} imply the continuity \ef{les2z} of
the solution $(\n,u)$, and the continuity of $\n{u}_y$ follows
from the equation~\ef{leq2} and that of $(\n,u)$. The proof of the
lemma is completed.
\enddemo
\par

Using Lemmas~\ref{lemma2}--\ref{lemma4} and a direct computation,
we can obtain the following result:
\oplem{}\label{lemma5} Let $T'_*>0$ be given in Lemma~\ref{lemma4}
and $(\n,u)$ be the solution of the IBVP
problem~\qef{l4a}--\qef{l4bc-b}. Then, under the assumptions of
Proposition~\ref{llocal},   $\n^{1+\al}{u}_y(t)\in C^0(\bar\Om)$
is uniformly bounded for any $t\in[0,T'_*]$
\be
 \lim_{y\to y_0}\n^{1+\al}u_y(y,t)=0,\quad
 \|\n^{1+\al}u_y\|_{L^\infty(0,T'_*;C^0(\bar\Om))}\le C(T'_*),\nnm
\ee
and $(\n,u)$ satisfies for $0\le s <t\le T'_*$ that
\bgr
 \|\n(t)-\n(s)\|_{L^2(\Om)}+ \|u(t)-u(s)\|_{L^2(\Om)}
 \le  C(T'_*)|t-s|,  \nnm
\\[2mm]
 \|\n^{1+\al}u_y(t)-\n^{1+\al}u_y(s)\|_{L^2(\Om)}
 \le  C(T'_*)|t-s|^{1/2}. \nnm
\egr
\cllem
\demo
The facts that for any $t\in[0,T_*]$, it holds that $
\n^{1+\al}{u}_y(t)\in C^0(\bar\Om)$ and $\lim_{y\to
y_0}\n^{1+\al}{u}_y(y,t)=0$ are due to the continuity of right hand
side terms of \ef{leq2a} and \ef{lasp2}. By \ef{les2} and
\ef{vacuum-point}, one can check easily that \be
 \|\n^{1+\al}u_y\|_{L^\infty(\Om\times[0,T'_*])}
 \le
 \|\n\|^{\al}_{L^\infty(\Om\times[0,T'_*])}
 \cdot\|\n{u}_y\|_{L^\infty(\Om\times[0,T'_*])}
 \le
  C(T'_*).  \nnm
\ee Making use of Eq.~\qef{l4a}, \ef{les1} and \ef{les2},  one can
obtain \bma & \|\n(t)-\n(s)\|_{L^2(\Om)}
 \le
  \int_s^t\|\n_t(\tau)\|_{L^2(\Om)}d\tau
 =\|\int_s^t\n^2u_y(\tau)\|_{L^2(\Om)}d\tau\nnm\\
 \le
&     C(t-s)\|\n\|_{L^\infty(\Om\times[0,T'_*])}
 \cdot\|\n{u}_y\|_{L^\infty(\Om\times[0,T'_*])}
 \le
  C(T'_*)|t-s|, \nnm
\ema
\be
  \|u(t)-u(s)\|_{L^2(\Om)}
 \le
  \|\int_s^tu_t(\tau)d\tau\|_{L^2(\Om)}
 \le
      C\int_s^t\|u_t\|_{L^2(\Om)}
 \le
  C(T'_*)|t-s|, \nnm
\ee
and
\bma
 &\|\n^{1+\al}u_y(t)-\n^{1+\al}u_y(s)\|_{L^2(\Om)}
\le
 \|\int_s^t(\n^{1+\al}u_y(t))_td\tau\|_{L^2(\Om)}\nnm\\
\le
 &C( \|\int_s^t\n_t\n^{\al}u_y(\tau)\|_{L^2(\Om)}
    +\|\int_s^t\n^{1+\al}u_{yt}(\tau)d\tau\|_{L^2(\Om)})\nnm\\
\le
 &C( \|\int_s^t\n^{2+\al}u^2_y(\tau)\|_{L^2(\Om)}
    +\|\int_s^t\n^{1+\al}u_{yt}(\tau)d\tau\|_{L^2(\Om)})\nnm\\
\le
 &C(T'_*)(|t-s|\|\n^{2+\al}u^2_y\|_{L^\infty(\Om\times[0,T'_*])}
    +|t-s|^{1/2}
    \|\n^{(1+\al)/2}u_{yt}(\tau)\|_{L^2(\Om\times[0,T'_*])})\nnm\\
\le
 &C(T'_*)|t-s|^{1/2}. \nnm
\ema The proof is completed.\enddemo

\bigskip

\underline{\it Proof of the Proposition~\ref{llocal}.} With the
help of Lemmas~\ref{lemma2}--\ref{lemma4}, we are ready to prove
the Proposition~\ref{llocal}.\par

{\it (1) Existence of weak solution  for short time.}
We only deal with the case for the Dirichlet boundary condition
and one point vacuum state in the initial data, the case of
periodic boundary and one point vacuum state in initial data can
be done in a similar way. Once the a-priori estimates are
established as in Lemmas~\ref{lemma2}--\ref{lemma5}, we are able
to prove the Proposition~\ref{llocal}. First of all, we construct
a sequence of approximate solutions by modifying the finite
difference scheme used in
\cite{Nishida86,LuoXinYang00,OkadaMatsusuMakino02}. Without the
loss of generality, we assume \be y_0=\frac12, \quad
\mbox{namely,}\quad \n_0(\mbox{$\frac12$})=0. \ee For any given
positive integer $N=2k+1$ with $k\ge 0$ an integer, let $h=1/N$.
Consider the system of $2N$ ordinary differential equations \be
\left\{\begin{aligned}
  &\frac{d}{dt}\n^h_{2n+1}
 +(\n^h_{2n+1})^2\frac{u^h_{2n+2}-u^h_{2n}}{h}=0,
\\[2mm]
  &\frac{d}{dt}u^h_{2n}
   +\frac{(\n^h_{2n+1})^\ga-(\n^h_{2n-1})^\ga}{h}\\
&\mbox{\quad} =\frac{1}{h}\{
    (\n^h_{2n+1})^{1+\al}(u^h_{2n+2}-u^h_{2n})/h
   -(\n^h_{2n-1})^{1+\al}(u^h_{2n}-u^h_{2n-2})/h \},
\end{aligned}\right.     \label{odea}
\ee
with the boundary condition and point vacuum
\be
 u^h_{0}(t)= u^h_{2N}(t)=0,\quad \n_{2k+1}(t)=0,
  \label{odeb}
\ee
and initial data
\be
 \begin{aligned}
  &\n^h_{2n+1}(0)=\n_0((2n+1)\mbox{$\frac{h}2$})
  &u^h_{2n}(0)=u_0((2n)\mbox{$\frac{h}2$})
\end{aligned}  \label{odec}
\ee
where $n=1,2,3,...,N$. Here, we also assume
\[
 \n^h_{2N+1}(t)=\n^h_{2N-1}(t),\ \
 \n^h_{1}(t)=\n^h_{-1}(t),\ \
 u^h_{2N+2}(t)=u^h_{-2}(t)=0
\]
which is consistent with the boundary condition and the fact that
density is continuous and non-zero at the left boundary.\par

By applying the idea as
in~\cite{LiuXinYang98,OkadaMatsusuMakino02} and the similar
arguments mentioned above, we can obtain the following uniform
(w.r.t. $h$) a-priori estimates about the solutions
$(\n_{2n+1},u_{2n})$ (here and below we omit the symbol $h$ for
simplicity) of \ef{odea}--\ef{odec} similar to
Lemmas~\ref{lemma2}--\ref{lemma5}. Details will be omitted (the
reader can refer to \cite{LiuXinYang98,OkadaMatsusuMakino02} for
similar arguments in details).
\par

\oplem{}\label{lem1}

Let $(\n_{2n+1},u_{2n})$ be the solution of
\qef{odea}--\qef{odec}. Then, it holds
\bgr
 \bln
 \sum_{n=0}^{N}&(\mbox{$\frac12$}u^2_{2n}(t) +\pi(\n_{2n+1}(t)))h
 +\int_0^t \sum_{n=1}^{N}\n^{1+\al}_{2n+1}
    (\mbox{$\frac{u_{2n+2}(s)-u_{2n}(s)}{h}$})^2hds\nnm\\
  &
  =
  \sum_{n=0}^{N}(\mbox{$\frac12$}u^2_{2n}(0)+\pi(\n_{2n+1}(0)))h,
\eln\\
 \sum_{n=0}^{N-1}\n^{-1}_{2n+1}(t)h
 =\sum_{n=0}^{N-1}\n^{-1}_{2n+1}(0)h.  \nnm
\egr \cllem It follows from Lemma~\ref{lem1} and  the standard
theory of ordinary differential equations that there exists a
global solution $(\n_{2n+1},u_{2n})$ to \ef{odea}--\ef{odec} for
any fixed positive $N$ and $h$. Furthermore, the following
properties hold:

\oplem{}\label{lem2}
Let $(\n_{2n+1},u_{2n})$ with $n=0,2,3,...,N-1$ with $N=2k+1$ be
the solution of \qef{odea}--\qef{odec}. Then, it holds for $n\ge
k+1$ that
\bgr
  \n^{\al}_{2n+1}(t)
 =
   \n^{\al}_{2n+1}(0)
  -\al\int_0^t\sum_{j=k+1}^{n}\frac{d}{dt}u_{2j}(s)hds
  -\al\int_0^t \n^\ga_{2n+1}(s)ds, \nnm
\\
 \n^{1+\al}_{2n+1}(t)\frac{u_{2n+2}(t)-u_{2n}(t)}{h} \nnm
 =
 \sum_{j=k+1}^{n}\frac{d}{dt}u_{2j}(t)h + \n^\ga_{2n+1},
\egr
and for $n\le k$ that
\bgr
  \n^{\al}_{2n-1}(t)
 =
   \n^{\al}_{2n-1}(0)
  -\al\int_0^t\sum_{j=n}^{k}\frac{d}{dt}u_{2j}(s)hds
  -\al\int_0^t \n^\ga_{2n-1}(s)ds, \nnm
\\
 \n^{1+\al}_{2n-1}(t)\frac{u_{2n}(t)-u_{2n-2}(t)}{h}
 =
  -\sum_{j=n}^{k}\frac{d}{dt}u_{2j}(t)h + \n^\ga_{2n-1}. \nnm
\egr
\cllem
\oplem{}\label{lem3}
Under the assumptions of Proposition~\ref{llocal}, there is a
short time $T'_*>0$ so that it holds for $n=0,1,2,...,N$ that
\bgr
  \n_{2n+1}(t)
 +
  \sum_{n=1}^{N}(
    \mbox{$\frac{\n^{\al}_{2n+1}(t)-\n^{\al}_{2n-1}(t)}
         {h}$})^2h
 +|\n_{2n+1}(t)\mbox{$\frac{u_{2n+2}(t)-u_{2n}(t)}{h}$}|
 \le  C(T'_*),                    \nnm
\\[2mm]
  a_-|\mbox{$\frac12$}(2n+1)h-\mbox{$\frac12$}|^\beta
   \le \n_{2n}(t)
  \le a_+|\mbox{$\frac12$}(2n+1)h-\mbox{$\frac12$}|^\beta,  \nnm
\egr
and for $ m=1\mbox{ or } n$ that
\[
 \begin{aligned}
   \sum_{j=0}^{N}(\mbox{$\frac{d}{dt}$}u_{2j}(t))^{2m}h
  +\int_0^t\sum_{j=0}^{N}\n^{1+\al}_{2j+1}(s)
   (\mbox{$\frac{d}{dt}$}u_{2j}(s))^{2m-2}
   (\mbox{$\frac{\frac{d}{dt}u_{2n}(s)
            -\frac{d}{dt}u_{2n-2}(s)}
         {h}$})^2)h
  \le  C(T'_*).               \nnm
 \end{aligned}
\]
Here $C(T'_*)>0$ and  $a_\pm>0$ are constants.
\cllem
\oplem{}\label{lem4}
Under the assumptions of Proposition~\ref{llocal} and
Lemma~\ref{lem3}, it holds for $t\in[0,T'_*]$ that
\bgr
|\n^{1+\al}_{2n+1}(t)\mbox{$\frac{u_{2n+2}(t)-u_{2n}(t)}{h}$}|
 \le  C(T'_*),\nnm
\\
    |u_{2n}(t)|+
 \sum_{n=1}^{N} |u_{2n}(t)-u_{2n-2}(t)| \le  C(T'_*),\nnm
\\
 \sum_{n=1}^{N} | (\n_{2n+1})^{1+\al}(u_{2n+2}-u_{2n})/h
   -(\n_{2n-1})^{1+\al}(u_{2n}-u_{2n-2})/h|\le  C(T'_*),\nnm
\\
  \sum_{n=1}^{N} |\n_{2n}(t) -\n_{2n}(s)|^2h
  +\sum_{n=1}^{N}|u_{2n-1}(t)-u_{2n-3}|^2h
  \le  C(T'_*)|t-s|^2,\nnm
\\
 \sum_{n=1}^{N} |\mbox{$(\n_{2n})^{1+\al}(t)
   \frac{u_{2n+1}(t)-u_{2n-1}(t)}{h}$}
   -(\n_{2n-2})^{1+\al}(s)
     \mbox{$\frac{u_{2n-1}(s)-u_{2n-3}(s)}{h}$}|^2h
   \le  C(T'_*)|t-s|.\nnm
\egr
\cllem

With the help of Lemmas~\ref{lem1}--\ref{lem4}, we can define the
sequence of approximate solutions $(\n_{h},u_{h})$ on the domain
$\Om\times[0,T'_*]$ as \be
\begin{cases}
 \n_h(y,t)=\n_{2n+1}(t),\\[2mm]
 u_h(y,t)=\frac1{h}\mbox{ $
            [(y-(n-\frac12)h)u_{2n+1}(t)
           +((n+\frac12)h-y)u_{2n-1}(t)]$}
\end{cases}      \nnm
\ee for $y\in((2n)\frac12h,(2n+2)\frac12h)$. It can be verified
that the following properties hold for the approximate solutions
\be
 \partial_yu_h(y,t)= \frac{u_{2n+1}(t)-u_{2n-1}(t)}{h} \nnm
\ee
and
\be
\begin{cases}
  a_-|y-y_0|^\beta\le \n_h(y,t) \le  a_+|y-y_0|^\beta,
\\[2mm]
 |u_h(y,t)|\le  C(T'_*),
\quad
 |\n_h(y,t)\partial_xu_h(y,t)|\le  C(T'_*), \nnm
\end{cases}
\ee for $(y,t)\in\Om\times[0,T'_*]$. Then, the existence of weak
solution for short time $t \in[0,t'_*]$ follows from  Helly's
theorem, the diagonal process together with Lebesgue's theorem,
and the a-priori estimates (see, for instance,
\cite{LuoXinYang00}).  The details are omitted.
\par\smallskip

{\it (2) Uniqueness of weak solution for short time}. We prove the
uniqueness of weak solutions for
$\al\in(\mbox{$\frac12$},\gamma)$. Without loss of generality,
only the case for the Dirichlet boundary condition \ef{l4bc-b}
will be studied. Let $(\n_1,u_1)$ and $(\n_2,u_2)$ be two weak
solutions of Eq.~\qef{l4a}--\qef{l4bc-b} satisfying
Lemma~\ref{lemma2}--\ref{lemma5}. Denote \be n=\n_1-\n_2,\quad
\psi=u_1-u_2,\quad   \nnm (x,t)\in(0,1)\times[0,T]. \ee Obviously,
the new unknown $(n,\psi)$ with $\psi(0,t)=\psi(1,t)=n(y_0,t)=0$
satisfies \bgr
 \left(\frac{n}{\n_1\n_2}\right)_t
 +  \psi_{y} = 0, \quad y\neq y_0,\label{l6a}
\\
 \psi_t
  + (\n_1^\ga-\n_2^\ga)_y
  - (\n_1^{1+\al}\psi_{y})_{y}
  - ((\n_1^{1+\al}-\n_2^{1+\al})u_{2y})_y=0,  \label{l6b}
 \egr
for $(y,t)\in(0,1)\times(0,T_*]$ with zero initial data \be
(n(y,0),\psi(y,0))=(0,0),\quad y\in (0,1).  \label{l6c} \ee Take
inner product between $\n_1^\al\n_2^{-1}n$ and \ef{l6a} over
$[0,y_0)\cup(y_0,1]$ to obtain
 \bma
 \frac{d}{dt}& \oint \n_1^{-1+\al}\n_2^{-2}n^2 dy
 =
   -(1+\al)\oint u_{1y}\n_1^{\al}\n_2^{-2}n^2dy
   -2\oint \n_1^{\al}\n_2^{-1}n\psi_{y}dy \nnm\\
 \le
  & (C\|\n_1u_{1y}\|_{L^\infty}+4)\oint \n_1^{-1+\al}\n_2^{-2}n^2dy
   +\frac14\oint \n_1^{1+\al}\psi^2_{y}dy\nnm\\
 \le
  &  C(T_*)\oint \n_1^{-1+\al}\n_2^{-2}n^2dy
    +\frac14\oint \n_1^{1+\al}\psi^2_{y}dy,     \label{lun1}
\ema where $\oint f dx =:\int_0^{y_0}fdx +\int_{y_0}^1fdx$. Taking
inner product between $\psi$ and \ef{l6b} over $\Om$, and noting
that $\psi(0,t)=\psi(1,t)=n(y_0,t)=0$ and
$\n_i^{1+\al}u_i(y_0,t)=0,i=1,2,$ we have \bma
\frac12\frac{d}{dt}&\|\psi(t)\|^2
   +\oint\n_1^{1+\al}\psi_{y}^2dy
 =
   \oint(\n_1^\ga-\n_2^\ga)\psi_{y}dy
   -\oint(\n_1^{1+\al}-\n_2^{1+\al})u_{2y}\psi_{y}dy
\nnm\\
 =
  &\oint\frac{\n_1^\ga-\n_2^\ga}{\n_1-\n_2}n\psi_{y}dy
   -\oint\frac{\n_1^{1+\al}-\n_2^{1+\al}}{\n_1-\n_2}nu_{2y}\psi_{y}dy
\nnm\\
\le
  &\, C\oint \left(\frac{\n_1^\ga-\n_2^\ga}{\n_1-\n_2}\right)^2
           \n_1^{-1-\al} n^2 dy
   + \frac18\oint \n_1^{1+\al}\psi^2_{y} dy\nnm\\
  &\,+4\oint\left(\frac{\n_1^{1+\al}-\n_2^{1+\al}}{\n_1-\n_2}\right)^2
          \n_1^{-1-\al}n^2u_{2y}^2 dy
   +\frac18\oint \n_1^{1+\al}\psi^2_{y}dy
\nnm\\
\le
  &\, \frac14\int \n_1^{1+\al}\psi^2_{y}dy
    + C\oint \n_1^{-1+\al}\n_2^{-2} n^2 dy.\label{lun3}
\ema Summing up the two differential inequalities \ef{lun1} and
\ef{lun3} leads to \be \bln
 &\frac{d}{dt}\|\psi(t)\|^2
    + \frac{d}{dt}\oint\n_1^{-1+\al}\n_2^{-2} n(t)dx
   +\oint\n_1^{1+\al}\psi_{y}^2dy \\
 \le
 &  C(T_*)\oint \n_1^{-1+\al}\n_2^{-2} n^2 dy
\eln \ee which, together with the initial data \ef{l6c} and the
Gr\"onwall's lemma, gives rise to \be \|\psi_y(t)\|^2
    + \oint\n_1^{-1+\al}\n_2^{-2} n^2(x,t)dx\equiv 0,\quad t\in[0,T].
\ee This together with $\n_1(y_0,t)=\n_2(y_0,t)=0$ implies the
uniqueness of weak solution \be (\n_1,u_1)\equiv (\n_2,u_2). \ee
The proof of Proposition~\ref{llocal} is completed.\enddemo
\par
\bigskip

\par\bigskip

\underline{\it Proof of Proposition~\ref{llocal-z}.}\ \
With the help of Propositions~\ref{llocal}--\ref{llocal-zz} and
the inverse transformation from the Lagrangian coordinates to the
Eulerian coordinates, we are able to prove the
Proposition~\ref{llocal-z}. Since the case of periodic boundary
conditions with one point vacuum state in the initial data can be
dealt with in a similar framework, we only deal with the IBVP
problem for the compressible Navier-Stokes
equations~\ef{l0a}--\ef{l0c} for the case of the Dirichlet
boundary condition \ef{d6} with one point vacuum in the regular
initial data \ef{l0h}--\ef{l0e}.

First, one can check easily that the IBVP
problem~\ef{l0a}--\ef{l0c}, \ef{d6}, and \ef{l0h}--\ef{l0e} in the
Eulerian coordinates is equivalent to the corresponding IBVP
problem in the Lagrangian coordinates for the
equations~\ef{l4a}--\ef{l4c} with the Dirichlet boundary
conditions \ef{l4bc-b} and \ef{linit2} through the coordinate
transformation \be
 y=\int_{0}^{x}\n(z,t)dz,
 \quad x\in(0,1), \ \ t\ge 0.  \label{tran-fa}
\ee Note that $y\in [0,1]$ due to the conservation of mass. More
importantly, the case of one point vacuum state in the initial
data \ef{l0h}--\ef{l0e} is reformulated into the corresponding
case \ef{linit2} with $y_0=\int_{0}^{x_0}\n_0(z)dz$ and
$\beta=\frac{\sigma}{1+\sigma}$. It is easy to verify that all the
assumptions of Proposition~\ref{llocal} are satisfied. Then
Proposition~\ref{llocal} gives the existence and uniqueness of
weak solution $(\tilde{\n},\tilde{u})$
satisfying~\ef{linit3a}--\ef{l4bc-az} in the Lagrangian coordinate
to the compressible Navier-Stokes equations~\qef{l4a}--\qef{l4c}
with the Dirichlet boundary conditions \qef{l4bc-b} and one point
vacuum state in the initial data. In particular,
\bgr
 \tilde{\n}, \tilde{u} \in C^0([0,T'_*];C^0(\bar\Om)),\quad
 \tilde{\n}{\tilde{u}}_y(t)
        \in L^\infty(0,T'_*;C^0(\bar\Om)),  \label{cor1}\\
 \tilde{\n}(y_0,t)=0,\quad  \n(y,t)>0,\ y\neq y_0.  \label{cor2}
\egr These in turn imply, in terms of the inverse coordinate
transformation of \ef{tran-fa}, i.e.,
\[
x =\int_0^y\tilde{\n}^{-1}(z,t)dz,\quad y\in[0,1],
\]
the existence and uniqueness of weak solution
$(\tilde{\n},\tilde{u})$ on the domain
$(x,t)\in[0,1]\times[0,T'_*]$ to the compressible Navier-Stokes
equations~\qef{l0a}--\qef{l0c} with the Dirichlet boundary
conditions \qef{d6} and initial one point vacuum state
\ef{l0h}--\ef{l0e}. Moreover, one can check that there exists one
particle path $x=X_0(t)$ defined by \be \
\dot{X}_0(t)=\tilde{u}(X_0(t),t),\ t> 0,\quad X_0(0)=x_0    \no
\ee satisfying
 \[
   y_0\equiv \int_{0}^{X_0(t)}\tilde{\n}(z,t)dz=
\int_{0}^{x_0}\n_0(z)dz
 \]
due to the fact
\[
 \frac{d}{dt}y(X_0(t),t)
 = \frac{d}{dt}\int_{0}^{X_0(t)}\tilde{\n}(z,t)dz=0,\quad t\ge 0.
\]
This, together with \ef{cor2},  gives rise to \be
 \tilde{\n}(X_0(t),t)=0,\quad\mbox{and}\quad
 \tilde{\n}(x,t)>0, \ x\neq X_0(t),\ t\ge 0.  \no
\ee Moreover, it is easy to verify that the solution
$(\tilde{\n},\tilde{u})$  satisfies all the properties
\ef{density}--\ef{short-b}, and particularly \be
\begin{cases}
  \tilde{\n}(x,T_*)\ge 0 \mbox{ on }\Om,\quad
    \tilde{m}(.,T_*)= \tilde{\n}\tilde{u}(.,T_*)=0,\
     \mbox{ on }\{x\in \Om\,|\,\n_0(x)=0\},
\\[2mm]
   \tilde{\n}(.,T_*)\in L^1(\Om)\cap L^\ga(\Om),\quad
   (\tilde{\n}^{\al-1/2}(.,T_*))_x\in L^2(\Om),
\\[2mm]
  \mbox{$\frac{|\tilde{m}(.,T_*)|^2}{\tilde{\n}(x,T_*)}$}
  +\mbox{$\frac{|\tilde{m}(.,T_*)|^{2+\nu}}
               {\tilde{\n}^{1+\nu}(x,T_*)}$}  \in L^{1}(\Om).
\end{cases}                   \nnm
\ee The case of periodic boundary can be treated in a similar way,
we omit the details. The proof of Proposition~\ref{llocal-z} is
completed.
\enddemo

\section{Vanishing of vacuum states and blow-up phenomena}
\setcounter{equation}{0}\label{Vanishing}
We shall prove that for any global  entropy weak solution $(\n,u)$
to the IBVP for compressible Navier-Stokes equations
\qef{l0a}-\qef{l0c} together with \qef{d6} or \qef{d7}, any
possible vacuum state vanishes in finite time and the velocity (if
definable and regular enough) blows up in finite time if vacuum
state appears, for example, the density contains vacuum initially
as in Theorem 2.2. The weak solution becomes a strong one after
the vanishing of vacuum states and tends time-asymptotically to
the non-vacuum equilibrium state exponentially.
\par

\subsection{Vanishing of vacuum states in finite time}
\oppro {}\lb{tt2} Let \ef{vanishing-vacuum} hold. For any global
entropy weak solution $(\n,u)$ to  the IBVP for compressible
Navier-Stokes equations \qef{l0a}-\qef{l0b}  with initial data
-\qef{l0c} and boundary condition \qef{d6} or \qef{d7} in the
sense of Definition~\ref{definition1}, there exists a time $T_0>0$
such that \be \inf\limits_{x\in\ol\Om}\n(x,t)>0,\mbox{ for all
}\quad t\ge T_0.\lb{f200} \ee \clpro

\demo To prove \ef{f200}, we will employ  an idea which has been
used in \cite{LiXin04} (see also \cite{li04,LiXinHuang05}) to show
the blow up behavior of  both the global strong solutions to the
IBVP for \ef{l0a}-\ef{l0b} with the constant viscosity and the
global strong solutions to the Stokes approximation equations, with
initial data containing vacuum states.

Let $T\in (0,\infty)$ be fixed. In this subsection, $C$ denotes
some generic positive constant independent of $T$. First, it is
noted that the total mass is conserved the total for any $t\in
(0,T]$
 \be
  \int_\Om \n (x,t)dx=\int_\Om \n_0 (x )dx. \lb{d101}
\ee Based on the entropy inequality \ef{d4}, it can be deduced
from \ef{d101} and \ef{d100} that for a constant
$b\ge\max\{\al+\ga-1 ,2\al+1 ,1\}$,
 \be
 \sup\limits_{0\le t\le T}
  \left(\|\n\|_{L^\infty}
        +\left\|\left(\n^b\right)_x\right\|_{L^2}\right)
 +\int_0^T\left\|\left(\n^b\right)_x\right\|_{L^2}^2dt
\le C.\lb{d12} \ee I
t will be shown below that
 \bn
 g(t)\triangleq\left\|\left(
  \rho^b-\ol{\rho^b}\right)(\cdot,t)\right\|^4_{L^{4}(\Om)}
\rightarrow 0\mbox{ as }t\rr \infty,\lb{d11}
 \en
where
\[
 \ol{\rho^b}(t)=\frac{1}{|\Om|}\int_{\Om}\rho^b(x,t) dx.
\]
Now, we assume that \ef{d11} holds,  and continue the proof of
Proposition \ref{tt2}. In fact, the inequality \ef{d12} and the
Poinc\'are-Sobolev inequality imply that
\bma
  \left\|\left(\rho^b-\ol{\rho^b}\right)(\cdot,t)\right\|_{C(\ol\Om)}
 &
 \le
  C\left\|
     \left(\rho^b-\ol{\rho^b}\right)(\cdot,t)
   \right\|^{2/3}_{L^{4}(\Om)}
   \|(\rho^b)_x(\cdot,t) \|^{1/3}_{L^2}
 \no&
 \le
  C\left\|
    \left(\rho^b-\ol{\rho^b}\right)(\cdot,t)
    \right\|^{2/3}_{L^{4}(\Om)}
  \rightarrow 0,\mbox{ as }t\rr\infty.\lb{d13}
 \ema
This suffices to finish the proof of Proposition \ref{tt2} due to
the following simple fact
\[
 \ol{\n^b}(t)\ge {\ol\n}^{\,b}(t)
  \equiv\ol{\n_0}^{\,b}=1,\mbox{ for any }t\ge 0.
\]

It remains to prove \ef{d11}.   First, it follows directly from
\ef{d12} and the Poinc\'are-Sobolev inequality
that
 \bn
 \int_0^Tg(t)dt&\le& C\sup\limits_{0\le t\le T}\left
\|\rho^b-\ol{\rho^b}\right\|_{L^\infty}^{2 }
 \int_0^T\|(\rho^b)_x\|_{L^2}^2dt  \le  C.\lb{d14}\en

Next, we prove that
 \bn
 \int_0^T|g'(t)|dt \le C.\lb{a72}
 \en
Note that \ef{d1} as well as the boundary condition \ef{d6} or
\ef{d7} imply
 that
  \bn
  g'(t)
   &= &
    4b\left\langle
       \left(\rho^b-\ol{\rho^b}\right)^{3}\rho^{b-1},\rho_t
     \right\rangle_{H^1\times H^{-1}}
   -4\left(\ol{\rho^b}\right)_t\int_\Om
    \left(\rho^b-\ol{\rho^b}\right)^{3}dx \no& =&
   -4b \int_\Om
      \left(
        \left(\rho^b-\ol{\rho^b}\right)^{3}\rho^{b-1}
       \right)_x\sqrt{\rho}\sqrt{\rho}udx
   -4\left(\ol{\rho^b}\right)_t\int_\Om
     \left(\rho^b-\ol{\rho^b}\right)^{3}dx
   \no&\triangleq&I_1+I_2.\lb{a55}
 \en
It follows from \ef{d12} and \ef{d4} that
 \bn
 \int_0^T|I_1|dt
&\le&C\int_0^T\left|\int_\Om
 \left(\rho^b-\ol{\rho^b}\right)^{2   }(\rho^b)_x\rho^{b-1/2}
 \sqrt{\rho} udx\right|dt\no&& +C\int_0^T\left|\int_\Om
 \left(\rho^b-\ol{\rho^b}\right)^{3}
       (\rho^{b-1})_x \sqrt{\rho}\sqrt{\rho}
 udx\right|dt\no&\le&C\int_0^T\left|\int_\Om
 \left(\rho^b-\ol{\rho^b}\right)^{2 }(\rho^{\al-1/2})_x
 \n^{b-\al }\sqrt{\rho}udx\right|dt
 \no&\le& C\int_0^T\left\|(\rho^{\al-1/2})_x\right\|_{L^2}\|
 \sqrt{\rho}u\|_{L^2}
 \left\|\rho^b-\ol{\rho^b}\right\|_{L^\infty}^{2 }dt
 \no&\le& C\int_0^T\|(\rho^b)_x\|_{L^2}^2 dt\no&\le&C.\lb{d16}
\en
 The uniform entropy estimate \ef{d4}, together with \ef{d1}, gives
that
 \bn
  \sup\limits_{0\le t\le T}\left|\frac{d}{dt}\ol{ \rho^b}(t)\right|
&=&
  b\sup\limits_{0\le t\le T}
  \left|\left\langle \rho^{b-1},\rho_t \right\rangle\right|\no
&=&
  b\sup\limits_{0\le t\le T}
  \left|\int_\Om\left(\rho^{b-1}\right)_x
                \sqrt{\rho}\sqrt{\rho} udx\right|\no
&\le&
 C\sup\limits_{0\le t\le T}
 \left|\int_\Om\left(\rho^{\al-1/2}\right)_x
                \rho^{b-\al}\sqrt{\rho}u dx\right|\no
&\le&
 C\sup\limits_{0\le t\le T}
 \left(\left\|\left(\rho^{\al-1/2}\right)_x
        \right\|_{L^2}\|\sqrt{\rho}u\|_{L^2}
 \|\rho\|_{L^\infty}^{b-\al}\right)\no
 &\le&
  C.             \lb{a73}
\en
 This together with Poincar\'e inequality and \ef{d12} yields
\bn
 \int_0^T|I_2|dt
\le
 C\int_0^T\left\|\rho^b-\ol{\rho^b}\right\|_{L^\infty}^{3}dt
\le
  C\int_0^T\|(\rho^b)_x\|_{L^2}^2dt
\le
  C.\lb{d17}
\en The estimate \ef{a72} thus follows directly from
\ef{a55}-\ef{d17}. Hence the desired estimate \ef{d11} follows
from \ef{d14} and \ef{a72}.  The proof of Proposition~\ref{tt2} is
completed.
\enddemo

\subsection{Regularity and asymptotics of weak solutions for large time}
\setcounter{equation}{0}

It is usually difficult to get information about the velocity field
for the global  entropy weak solution $(\n,\sqrt{\n}u)$, in the
sense of Definition~\ref{definition1} to the IBVP   for the
Compressible Navier-Stokes equations \qef{l0a}-\qef{l0b} with
initial data \qef{l0c} and boundary values \ef{d6} or \ef{d7} in the
appearance of vacuum states. After vacuum states vanish, however, it
will be shown that the velocity field $u$ can be defined with enough
regularity and the nonlinear diffusion term is represented in terms
of the velocity $u$ and the density $\n$. The momentum equation
becomes a uniform parabolic equation, and the weak solution
$(\n,\sqrt{\n}u)=(\n,\sqrt{\n}\cdot u)$ becomes a strong solution.

Proposition \ref{tt2} implies that there is a time $T_0>0$ after
which the density of the global  entropy weak solution $(\n,u)$ to
the IBVP problem for \qef{l0a}-\qef{l0c} together with   \ef{d6}
or \ef{d7} is strictly positive  and  $(\n,u)$  satisfies the
finite entropy estimate~\ef{d4}. Consider the IBVP problem
\ef{l0a}--\ef{l0b}
  again for time $t\ge T_0$ with data
given at time $t=T_0$ by \be \bln
 \n(x,T_{0})=\lim_{t\to T_{0 } \mbox{}}\n(x,t),\quad
 {u}(x,T_{0})=\lim_{t\to T_{0 }}
  \mbox{$\frac{\sqrt{\n}u(x,t)}{\sqrt{\n(x,t)}}$},
   \label{l0ce}
\eln \ee and note here that away from vacuum the Dirichlet boundary
condition \ef{d6} reduces to \be
 u(0,t)=u(1,t)=0,\quad t\ge T_0. \label{d6b}
\ee We then have the regularity property of the solution for the
compressible Navier-Stokes equations~\qef{l0a}--\qef{l0c} with the
Dirichlet boundary condition \ef{d6} or the periodic boundary
condition~\ef{d7} for large time.
\oppro{} \label{lstrong}
Under the assumptions of Theorem~\ref{entropy_solution}, let
$(\n,\sqrt{\n}u)$ be the global  entropy weak solution to  the
IBVP for the compressible Navier-Stokes
equations~\qef{l0a}-\qef{l0b} with initial data \qef{l0c} and
boundary value \qef{d6} or \qef{d7} in the sense of
Definition~\ref{definition1}. Let $T_0>0$ so that the global weak
solution $(\n,u)$  satisfies for two positive constants $\n_\pm$
that
 \be
  0<\n_-\le \n(x,t)\le \n_+,
 \quad \forall\ (x,t)\in\bar\Om\times[T_0,\infty). \label{lstrong-2a}
\ee Then, $(\n,\sqrt{\n}u)=(\n,\sqrt{\n}\cdot u)$ is the unique
strong solution to the IBVP for the compressible Navier--Stokes
equations~\qef{l0a}--\qef{l0b} and \qef{l0ce} with the boundary
condition~\qef{d6}\footnote{Note here that away from vacuum the
Dirichlet boundary condition \ef{d6} reduces to the usual one
$u(0,t)=u(1,t)=0,\ t\ge T_0$.} or \qef{d7} for $t\ge T_0$.
Moreover, the regularity \qef{vvs3b} and the long time behavior
\qef{vvs4} hold. \clpro

\demo We only prove Proposition~\ref{lstrong} for the Dirichlet
case below, the periodic case can be treated similarly.\par

 \underline{\it Step 1. Regularity}.\
 It follows easily from proposition
\ref{tt2} that there exist  some $T_0$ and a constant $\rho_->0$
such that for all $t\ge T_0,$
 \be
\inf\limits_{x\in\ol\Om}\n(x,t)\ge \rho_->0 \lb{f1}
 \ee
which, together with \ef{d4}, implies \be
 \n  \in L^\infty (T_0,T; H^1(\Om)) \label{f1z}
 \ee
for any $T>T_0$. By the continuity of $\n\in
C(\bar\Om\times[0,\infty))$ there exists some $\si>0$ small enough
such that for any $t\ge T_0-\si,$
\bnn
 \inf\limits_{x\in\ol\Om}\n(x,t)\ge \frac{\rho_-}{2}>0.
 \enn
This implies that one can define the velocity $u$ for any global
 entropy weak solution in the sense of Definition~\ref{definition1}
after the vanishing of vacuum states by
 \[
 u=: \dis{\frac{\sqrt{\n}u}{\sqrt{\n}}}, \quad \mbox{for}\quad
t\ge T_0-\si.
\]
It then follows from the definition and \ef{s1} that
 \be
 u\in  L^\infty(T_0-\si,T;L^2(\Om)) \lb{f54}
\ee for any $T>T_0$.  Noting that \ef{f50} implies that for any
$\vp(x)\in C_0^\infty(\Om), \psi(t)\in C_0^\infty (T_0-\si,T),$
\bnn \lefteqn{
 \int_{T_0-\si}^T\psi(t)\int_\Om \Lambda\n^{-\al}\vp dx dt}\no
&&=
 -\int_{T_0-\si}^T
  \psi(t)\int_\Om\rho^{\al-1/2}\sqrt{\n}u(\n^{-\al}\vp)_xdxdt\no
&&\aaa
 -\frac{2\al}{2\al-1}\int_{T_0-\si}^T
  \psi(t)\int_\Om(\rho^{\al-1/2})_x\sqrt{\n}u \n^{-\al}\vp  dxdt\no
&&
 =-\int_{T_0-\si}^T\psi(t)\int_\Om\sqrt{\rho} u \n^{-1/2}\vp_xdxdt\no
&&
 =-\int_{T_0-\si}^T\psi(t)\int_\Om u\vp _x dx dt,
\enn
we can define the spatial  derivative of velocity and, together
with \ef{d4}, its regularity as
 \be
 u_x=\frac{\Lambda}{\n^\al} \in L^2(\Om\times(T_0-\si,T)).\lb{f3}
 \ee
In terms of \ef{f54}, \ef{f3} and \ef{d1} we are also able to
justify the Dirichlet boundary condition \ef{d6} for the velocity
$u$
 \bn
 u(0,t)=u(1,t)=0, \mbox{ for any }t\ge T_0-\si.\lb{f4}
\en Thus, \ef{f54}, \ef{f3} and \ef{f4} show \bn
 u\in L^2(T_0-\si,T;H^1_0(\Om))
    \cap L^\infty(T_0-\si,T;L^2(\Om)) \lb{f5}
 \en
for the case of the Dirichlet boundary conditions. Note here that
$u\in L^2(T_0-\si,T;H^1_{\rm per}(\Om))
    \cap L^\infty(T_0-\si,T;L^2(\Om))$ in the case of
periodic boundary conditions. We thus obtain from \ef{d1} and
\ef{d2} that the solution $(\n,u)$ satisfies \be
 \n_t+(\n u)_x=0\quad \mbox{ a.e. in }\quad
 \Om\times(T_0-\si,T),\lb{f6}
 \ee
and
 \be
  \int_{T_0-\si}^T\int_\Om \n u\vp_tdxdt
 +\int_{T_0-\si}^T\int_\Om
 \left(\n u^2-\n^\al u_x+\n^\ga\right) \vp_xdxdt=0  \lb{f7}
\ee for any $ {\vp}(x,t)\in C_0^\infty(\Om\times(T_0-\si,T))$ for
the Dirichlet case. The Eq.~\ef{f7} can be re-written in terms of
\ef{f5}, \ef{f6} and \ef{d4} as follows \bn
    \int_{T_0-\si}^T\int_\Om (u\vp_t   - \n^{\al-1}u_x\vp_x
  + \left( \n^{\al-2}\n_x-u \right)u_x\vp) dxdt
 =\int_{T_0-\si}^T\int_\Om\ga\n^{\ga-2}\n_x\vp\, dxdt\lb{f55}
 \en
for any $ {\vp}(x,t)\in C_0^\infty(\Om\times(T_0-\si,T))$ for the
Dirichlet case.

Noticing that $  \n^{\al-2}\n_x-u \in
L^\infty(T_0-\si,T;L^2(\Om))$ due to \ef{s1}, \ef{f1} and
\ef{f54}, and using standard regularity results for linear
parabolic equations (see \cite{la1}),  we get that \bn u\in
L^2(T_0,T ;H^2(\Om)) )\cap H^1(T_0,T;L^2(\Om)) \lb{f56}
 \en
 for the Dirichlet case.
It is noted here that $u\in L^2(T_0,T;H^2_{\rm per}(\Om))
    \cap H^1(T_0,T;L^2(\Om))$ for
the periodic case.

\underline{\it Step 2. Uniqueness}.\  We shall show that if there
exists another solution $(\eta,v)$ to the compressible
Navier-Stokes equations \ef{l0a}-\ef{l0b} with the following
initial data and Dirichlet boundary conditions
 \bn
 \begin{cases}
  (\eta,v)(x,T_0)=(\n,u)(x,T_0),\\
  v(0,t)=v(1,t)=0,
  \end{cases} \lb{f25}
\en
such that
\be
  \begin{cases}
  \n_-\le\eta\in L^\infty(T_0,T;H^1(\Om)),\\
   v\in L^\infty(T_0, T ;H ^1(\Om))\cap
        L^2(T_0,T;H^2(\Om))\cap H^1(T_0,T;L^2(\Om)),
  \end{cases}\lb{f13}
 \ee
then
 \be
    \n=\eta,\quad u=v\quad  \mbox{ a.e. in }\quad
          \Om\times (T_0,T ). \lb{f12}
 \ee

In fact, it follows  from \ef{f6}, \ef{f7} and \ef{f56} that
 \be
  \frac{1}{2}\int_\Om \n u^2dx
 +\int_{T_0}^t\int_\Om \n^\al u_x^2dxds
 -\int_{T_0}^t\int_\Om \n^\ga  u_x dxds
 =
 \frac{1}{2}\int_\Om \n u^2(x,T_0)dx        \lb{f15}
 \ee
for all $t\in (T_0,T),$ while \ef{f7} and \ef{f13} imply that
 \bn
  \lefteqn{
    \int_\Om \n u vdx+\int_{T_0}^t\int_\Om \n^\al u_xv_xdxds
   -\int_{T_0}^t\int_\Om \n^\ga  v_x dxds}\no&&
  =
   \int_\Om\n u v(x,T_0)dx
     +\int_{T_0}^t\int_\Om \n u(v_t+uv_x)dxds    \lb{f8}
 \en
for all $t\in (T_0,T)$.  To estimate the second term on the right
hand side of \ef{f8}, we use the decomposition \bn
 \n v_t+\n uv_x
 =(\n-\eta)(v_t+vv_x)+\n(u-v)v_x
 +( \eta^\al v_x)_x-(\eta^\ga)_x.\lb{f16}
 \en
Multiplying \ef{f16} by $u$ and we integrating by parts give
 \bn
 \lefteqn{\int_{T_0}^t\int_\Om \n u(v_t+uv_x)dxds}\no&&
 =
  \int_{T_0}^t\int_\Om (\n-\eta)u(v_t+vv_x)dxds
 +\int_{T_0}^t\int_\Om \n u(u-v)v_xdxds\no&&\aaa
 -\int_{T_0}^t\int_\Om \eta^\al v_xu_xdxds
 +\int_{T_0}^t\int_\Om  \eta^\ga u_x dxds.\lb{f17}
\en
 Substituting \ef{f17} into \ef{f8} gives, for
a.e. $t\in (T_0,T),$ that
 \bn
  \lefteqn{\int_\Om \n u vdx+\int_{T_0}^t\int_\Om \n^\al u_xv_xdxds
   - \int_{T_0}^t\int_\Om(\n^\ga  v_x +\eta^\ga u_x)dxds}\no&&
 =
   \int_\Om \n u v(x,T_0)dx
  +\int_{T_0}^t\int_\Om (\n-\eta)u(v_t+vv_x)dxds\no&&
  \aaa
  +\int_{T_0}^t\int_\Om \n u(u-v)v_xdxds
  -\int_{T_0}^t\int_\Om\eta^\al v_xu_xdxds.\lb{f18}
 \en
Multiplying \ef{f16} by $v$ and integrating the result over
$\Om\times(T_0,t)$ lead to \bn  \frac12\int_\Om \n v^2dx
 &=&\frac12\int_\Om \n  v^2(x,T_0)dx +
 \int_{T_0}^t\int_\Om (\n
-\eta)v(v_t+vv_x)dxds\no&& +\int_{T_0}^t\int_\Om \n v(u-v)
v_xdxds-\int_{T_0}^t\int_\Om \eta^\al v_x^2dxds\no&&+
\int_{T_0}^t\int_\Om
 \eta^\ga v_x dxds.\lb{f19}\en

We obtain after adding up \ef{f15} and \ef{f19} and subtracting
\ef{f18} that  for all $t\in (T_0,T),$
\bn
 \lefteqn{\frac12\int_\Om \n (u-v)^2dx
  +\int_{T_0}^t\int_\Om \n^\al (u-v) ^2_xdxds}\no&&
 =
  \int_{T_0}^t\int_\Om (\n^\al-\eta^\al) (v-u)_xv_xdxds
 +\int_{T_0}^t\int_\Om (\n-\eta)(v_t+vv_x)(v-u) dxds\no&&\aaa
 -\int_{T_0}^t\int_\Om \n  (u-v)^2v_xdxds
 -\int_{T_0}^t\int_\Om
   \left(\n^\ga   -\eta^\ga \right)(v-u)_xdxds\no&&
\le
 C_\ve \int_{T_0}^t
    \left(\|\n^\al-\eta^\al\|_{L^2}^2\|v_x\|_{L^\infty}^2
  + \|\n -\eta \|_{L^2}^2\|v_t+vv_x\|_{L^2}^2
  + \|\n^\ga-\eta^\ga\|_{L^2}^2\right)ds\no&&\aaa
  + C\int_{T_0}^t\|v_x\|_{L^\infty}\int_\Om \n(u-v)^2dxds
  + C\ve\int_{T_0}^t \| u-v \|_{H^1}^2ds.\lb{f20}
 \en
Next, we estimate the term
 \[
 g(t)\triangleq\|\n -\eta
\|_{L^2}^2+\|\n^\al-\eta^\al\|_{L^2}^2+\|\n^\ga
-\eta^\ga\|_{L^2}^2.
 \]
For any $\beta>0,$ \ef{f13},   \ef{f56} and \ef{l0a} imply that
 \bn
  (\n^\beta-\eta^\beta)_t+v(\n^\beta-\eta^\beta)_x
 +(u-v)(\n^\beta)_x+\beta\n^\beta (u-v)_x
 +\beta(\n^\beta-\eta^\beta)v_x=0.\lb{f21}
 \en
One can derive from this that
 \bn
  \left(\|\n^\beta-\eta^\beta\|_{L^2}^2\right)_t
&\le&
 \|\n^\beta-\eta^\beta\|_{L^2}^2
  \left(C\|v_x\|_{L^\infty}+C_\ve\|\n_x\|^2_{L^2}
 +C_\ve\right)\no&&
 +\ve\|u-v\|_{L^\infty}^2+\ve\|(u-v)_x\|_{L^2}^2.\lb{f22}
 \en
Thus, since $g(T_0)=0,$ we obtain from \ef{f22} with
$\beta=1,\al,\gamma$ respectively that
\be
 g (t)
  \le
  \int_{T_0}^tg(s)\left(C\|v_x\|_{L^\infty}
 +C_\ve\|\n_x\|^2_{L^2}
 +C_\ve\right)ds
 +C\ve\int_{T_0}^t\| u-v \|_{H^1}^2ds.\lb{f23}
 \ee

Now \ef{f12} is a consequence from \ef{f20}, \ef{f23}, \ef{f13},
\ef{f1z} and \ef{f1}. The proof of large time convergence \ef{vvs4}
follows directly from the standard arguments (see \cite{stz1} for
instance) with the help of entropy inequality~\ef{d4}. The proof of
Proposition~\ref{lstrong} is completed.\enddemo

\subsection{Finite time blow-up}
\label{finite_time_blowup}
In this subsection, we shall prove the Theorem~\ref{blow-up} about
the finite time blow-up phenomena as an immediate consequence of
Theorem~\ref{entropy_solution}, Propositions~\ref{tt2} and
Proposition~\ref{lstrong}.

\underline{\it Proof of Theorem~\ref{blow-up}}.\  We will prove
\ef{blowup-1a} only. The proof of \ef{blowup-1g} is similar. If
\ef{blowup-1a} fails, then there exists a fixed constant $\eta>0,$
such that \be
  \int^{T_1+\eta}_{T_1}\|u_x\|_{L^\infty}ds<\infty.   \lb{ab3}
 \ee
For any $(x,t)\in \ol\Om\times(T_1,T_1+\eta]$,  the particle path
$x(s) = X(s;t,x)$ through $(x,t)$ is given by  \be
 \lb{ab1}
\begin{cases}
 \frac{\pa}{\pa s}X(s;t,x)=u(X(s;t,x),s),
   & T_1\le s<t\le T_1+\eta,\\
 X(t;t,x)=x, & T_1\le  t\le T_1+\eta,\  x\in \ol\Om,
\end{cases}
 \ee
which  is well-defined due to \qef{ab3} and \ef{vvs3b}.
Consequently, one obtains via a standard argument from the
transport equation \ef{l0a} \be
 \lb{ab2}
 \rho (x,t)
   =
    \rho(X(T_1;t,x),T_1)
    \exp\left\{-\int_{T_1}^t u_y(y,s)|_{y=X(s;t,x)}ds\right\}
 \ee
for any $(x,t)\in \ol\Om\times (T_1,T_1+\eta]$. On the other hand,
it follows from \ef{ab3} and \ef{ab1} that for any $x\in \ol\Om$
there exists a trajectory $x=x(t)\in \ol\Om$ for $t\in
[T_1,T_1+\eta]$ so that $X(T_1;t,x(t))=x$. In particular, there
exists a trajectory $x=x_1(t)\in \ol\Om$ for $t\in [T_1,T_1+\eta]$
so that $X(T_1;t,x(t))=x_1$ with $(x_1,T_1)$ determined by
\ef{blowup-1c}, namely, $\rho (x_1,T_1)=0$. Thus, due to
$(\ref{ab2}),$ we deduce from \qef{ab3} that
\[
\rho (x_1(t),t)\equiv 0\quad \mbox{  for all }\quad  t\in
(T_1,T_1+\eta],
\]
which contradicts \ef{blowup-1c}. Thus, the blowup phenomena
\ef{blowup-1a} happens.  The proof of the Theorem~\ref{blow-up} is
completed.\enddemo

\bigskip\vspace{1cm}
\noindent\textbf{Acknowledgements}  The authors would like to thank
the referee for informing them the recent published references
\cite{bd07,bdg07,mv07} and the helpful comments on the manuscript.
They also would like to thank Prof. Didier Bresch for his interests
in this work and helpful discussions about the BD entropy estimates
he and his collaborators introduced.

The main part of this research was done when H.Li and J. Li were
visiting the Institute of Mathematical Sciences (IMS) of The Chinese
University of Hong Kong.
 The financial supports form the IMS
and the hospitality of the staff at the IMS are appreciated greatly.
The research of H. Li is partially supported by Beijing Nova
Program, NNSFC No.10431060, the NCET support of the Ministry of
Education of China, the institute of Mathematics and
interdisciplinary Science at CNU, the Re Shi Bu Ke Ji Ze You
program, and Zheng Ge Ru foundation. The research of J. Li is
partially supported by the JSPS Research Fellowship for foreign
researchers, NNSFC No.10601059. The research of J. Li and Xin are
partially supported by Hong Kong RGC Earmarked Research Grants
CUHK4028/04 and 4040/02, and Zheng Ge Ru Foundation.

\end{document}